\magnification =\magstep1
\baselineskip =13pt

\overfullrule =0pt

\def\H{{ H^{^{^{\hskip -.17cm \bullet}}}}}
\def\cher{{H^{^{^{\hskip -.22cm \bullet\bullet}}}}}
\def\ind{{\lim\limits_{\longrightarrow}}}
\def\pro{{\lim\limits_{\longleftarrow}}}
\def\bd{{\big\downarrow}}

\centerline {\bf DOUBLE AFFINE HECKE ALGEBRAS}

\centerline {\bf  AND 2-DIMENSIONAL LOCAL FIELDS}

\vskip .7cm

\centerline {\bf M. Kapranov}

\vskip 1cm

The concept of an $n$-dimensional local field was introduced by
A.N. Parshin [Pa 1] with the aim of  generalizing
 the classical adelic formalism
to (absolutely) $n$-dimensional schemes. By definition, 
a 0-dimensional local field is just a finite field
and an $n$-dimensional local field, $n>0$, is a complete discrete valued
field whose residue field is $(n-1)$-dimensional local. Thus for $n=1$
we get the locally compact fields such as ${\bf Q}_p$, $F_q((t))$ and
for $n=2$ we get fields such as ${\bf Q}_p((t))$,
$F_q((t_1))((t_2))$ etc. 

\vskip .2cm

In representation theory, harmonic analysis on reductive
groups over 0- and 1-dimensional local fields leads, in particular,
to consideration of the finite and affine Hecke algebras
$H_q, \H_q$ associated to any finite root system $R$ and any $q\in {\bf C}^*$.
These algebras can be defined in several ways, one being by generators
and relations,  another as the convolution algebra, with respect to
the Haar measure, of functions on the group bi-invariant
with respect to an appropriate subgroup (i.e., as the algebra of double
cosets). Harmonic analysis on groups over 2-dimensional local fields
has not been developed, the main difficulty being 
 the infinite dimensionality (absense of local
compactness) of such fields. However,  the double affine Hecke algebra
$\cher_q$ recently defined by I. Cherednik [Ch] in terms of
generators and relations, looks like the third term in the
hierarchy starting from $H_q, \H_q$. The problem
``give a group-theoretic construction of the Cherednik algebra"
(i.e., realize it as some algebra of double cosets) was proposed by D. Kazhdan
a few years ago.

\vskip .2cm

The purpose of the present paper is to provide a solution to
this problem by developing beginnings of harmonic analysis on
reductive groups over 2-dimensional local fields. We consider
a simple algebraic group $G$ (over {\bf Z}),  a 2-dimensional local
field $K=k((t))$ of equal characteristic (so $k$ is 1-dimensional local)
and the canonical
central extension $\Gamma$ of $G(K)$ by $k^*$. For an appropriate subgroup
$\Delta_1\i\Gamma$ the fibers of the Hecke correspondences (3.1) are
locally compact spaces (affine spaces over $k$ of growing dimension)
which possess natural invariant measures, so one can
{\it formally} define the Hecke operators 
associated to double cosets, by integrating over these measures.
The main difficulty here is the noncompactness of the
domain of integration. It poses covergence problems, making
it unclear how to compose such operators or how to define their
action from some vector space to itself. More precisely, the
operators are well defined on the space ${\cal F}_0$ of functions
on $\Gamma/\Delta_1$ with certain proper support conditions
but their values lie in a bigger space $\cal F$. 

\vskip .2cm

The way around the difficulty that we take is to use the analytic continuation
with respect to the parameters of the principal series representations.
This is a version of the classical method of regularizing
divergent integrals by introducing complex powers of auxiliary
polynomials into the integrand, so that the regularized integral
comes out as a function (possibly meromorphic) of the complex exponents.
In our nonarchimedean case such  functions are rational,
 in appropriate variables. 
Then, we define the regularized Hecke algebra $H(\Gamma, \Delta_1)$
to consist
of certain linear combinations of the regularized Hecke operators
with coefficients being rational functions of the principal
series exponents. The rule to single out the admissible linear
combinations is that they should preserve the space ${\cal F}_0$.
The main result (Theorem 3.3.8) is that $H(\Gamma, \Delta_1)$
is isomorphic to the (slight modification of) the Cherednik algebra
associated to $G$.  

\vskip .2cm

These analytic difficulties seem unavoidable because we are
dealing with an infinite-dimensional group with no Haar measure. 
The approach presented here is not restricted to our particular
choice of the subgroup $\Delta_1$ but can be applied
to other choices as well, in particular, to arbitrarily ``deep"
congruence subgroups $\Delta'\i\Delta_1$ (with
$\Delta_1/\Delta'$ locally compact).  As in the familiar p-adic case,
the resulting algebras, while not possessing a nice independent
description, are nevertheless important for the general theory
of representations of infinite-dimensional groups such as $\Gamma$. 
They will be studied is another paper. 

Another general point which seems important for the future is
the essential role of iterated ind- and pro-objects
in the representation theory of groups over 2-dimensional local
fields. Thus, the ``spaces" of the principal series representations
are in fact pro-vector spaces, the group itself is not
a topological group at all but rather a group object in a certain
iterated pro-ind-category and so on. The philosophy that
topological concepts when applied to $n$-dimensional local fields,
$n\geq 2$, become inadequate and should be replaced by considering
ind/pro objects, was explicitly formulated by K. Kato [Kat]. 
In the present paper we use ind/pro-objects in a
systematic way. 

\vskip .2cm

In order to keep the paper short, it was organized as follows.
Section 1 contains the setup for the groups and homogeneous spaces
we consider, including the issues related to the central extension.
This material is mostly well known. Section 2 is a reminder on
the Cherednik algebras. In Section
3 we explain our approach to constructing the regularized Hecke
algebra and formulate the main theorem. We tried to
split the construction into several steps so that they are
easily generalizable to more complicated situations. 
Section 4 contains the proof of the main theorem by using the principal
series intertwiners and  a version of the Mellin transform.
Finally, the Appendix contains the general constructions
on ind/pro objects (such as function spaces etc.) necessary
for the main text. Its content is used throughout the paper. 

\vskip .2cm

Among the (already fairly numerous) works on the Cherednik
algebras two are most relevant for this paper. One is the
paper [GG] which gives a construction in terms of
equivariant K-theory, generalizing the results of [KL]. 
The other one is [GKV] where an ``elementary" description
was given by unraveling [KL] [GG] in terms of certain
residue conditions. This description is in fact used 
here in an essential way, as the conditions of [GKV] match
very precisely the singularity 
properties of the principal series intertwiners. 

\vskip .2cm

I would like to thank A. Braverman, J. Bernstein, S. Evens,
V. Ginzburg and A. Parshin for useful discussions and correspondence.
In particular, work with V. Ginzburg on a related project stimulated my
thinking about the questions addressed here. 
The work on this paper was partially supported by an NSF grant
and a large part of it was carried out during a visit to the IHES
in the summer of 1998.

\vskip 2cm

\centerline {\bf \S 1. P-adic loop groups and their homogeneous spaces.}

\vskip 1cm

 \noindent {\bf (1.1) Root systems.} We start by introducing notation,
to be used in the rest of the paper. Let $G$ be a split simple,
simply connected
algebraic group (over {\bf Z}), $T\i B\i G$ the fixed maximal
torus and Borel subgroup, $N=[B,B]$. We regard $G, B, N, T$ as
group schemes. 

Let $L={\rm Hom}({\bf G}_m, T)$ and $L^\vee = {\rm Hom}(T, {\bf G}_m)$
be the coweight and weight lattices of $G$; $R\i L^\vee$ the root
system, $R_{\rm sim}\i R_+\i R$ the sets of simple and
positive roots. For $\alpha\in R_+$ let $\alpha^\vee\in L$ be
the corresponding coroot. Let also $\theta\in R_+$ be the maximal root
and $\rho\in L^\vee$ be the half-sum of positive roots.
By $h^\vee$ we denote the dual Coxeter number of $G$. 

We denote by $\check T = {\rm Spec} \,\, {\bf C}[L]$ the complex torus
 dual to $T$.  Let $W$ be the Weyl group of $G$. It acts on $L$ and $L^\vee$.
Denote by
$$\Psi: L\to L^\vee, \quad \Psi(a) = {1\over h^\vee}
\sum_{\alpha\in R_+} (\alpha,  a)
\alpha
\leqno (1.1.1)$$
the minimal
 $W$-invariant integral scalar product on $L$, see
[Kac], \S 6.

Let $L_{\rm aff} = {\bf Z}\oplus L$, $L^\vee_{\rm aff} = {\bf Z}\oplus L^\vee$ be
the lattices of affine (co)weights of $G$. They are dual to each other.
Let
$$\widehat R = \{ (n,\alpha), \alpha\in R, n\in {\bf Z}\} \i L^\vee_{\rm aff},
\quad \widehat  R_+ = (\{0\}\times R_+) \,\, \cup \,\, ({\bf Z}_{>0}\times
 R)$$
be the system of affine roots of $G$ and the set of positive affine roots. The
set of simple affine roots will be denoted by
$$\widehat  R_{\rm sim} = (\{0\}\times R_{\rm sim}) \,\,\cup\,\, \{(1, -\theta)\}.$$
let $\widehat  W_{} = L \propto W$ be the affine Weyl group of $G$, generated by the
reflections $s_\alpha, \alpha\in\widehat  R_+$. 
We denote by $\leq$ the Bruhat order on $\widehat  W$. . Let
$l: \widehat  W\to {\bf Z}_+$ be the length function.
The group $\widehat  W$ acts {\bf Z}-linearly on $L^\vee_{\rm aff}$ by
$$w(m,b) = (m, w(b)), \,\, w\in W, \quad a(m,b) = (m+(a,b), b),\,\,a\in L.
\leqno (1.1.2)$$
We will also need the action on $L_{\rm aff}$ given by
$$w\circ (m,b) = (m, w(b)), w\in W, \quad
a\circ (m,b) = (m+ \Psi(a,b), b), a\in L. \leqno (1.1.3)$$
 The action (1.1.2)
preserves $\widehat R$, and we set
$$  D(w) = \{ \alpha\in \widehat  R_+: \,\, w(\alpha)\in\widehat R_-\}, \quad w\in\widehat  W.\leqno (1.1.4)$$
Then $| D(w)|=l(w)$. 
Let 
$$\widehat \rho = (-h^\vee, \rho) \in L^\vee_{\rm aff} \leqno (1.1.5)$$
be the standard substitute for the half-sum of positive affine roots,
see [PS] \S 14.3. It satisfies $(\widehat \rho, \alpha^\vee)=1$ for any
$\alpha\in\widehat  R_{\rm sim}$. We set
$$\widehat \delta(w) = 
\widehat \rho -w(\widehat \rho) = \sum_{\alpha\in   D(w)}\alpha, \quad
w\in\widehat  W. \leqno (1.1.6)$$
Let also
$$T_{\rm aff}({\bf C}) = {\rm Spec} \,\, {\bf C}[L^\vee_{\rm aff}], 
\quad \check T_{\rm aff} = {\rm Spec} \,\, {\bf C}[L_{\rm aff}]
\leqno (1.1.7)$$
to be the affine tori corresponding to $T({\bf C})$, $\check T$. 

\vskip .3cm

\noindent {\bf (1.2) Groups and homogeneous spaces.} 
 Let $k$ be a complete discrete valued field with residue field
$F_q$ and $K=k((t))$. So 
 $K$ is a complete discrete valued field with residue field $k$, thus
 a 2-dimensional local field in the sense of [Pa1]. We
use the notations ${\cal O}_K = k[[t]], {\bf m}_K = tk[[t]]$ for the ring of
integers and maximal ideal of $K$  and ${\cal O}_k, {\bf m}_k$
for the analogous subrings in $k$. 
Denote $\pi_K: {\cal O}_K\to k$, $\pi_k: {\cal O}_k\to F_q$ the projections
 and set
${\cal O}' = \pi_K^{-1}({\cal O}_k)\i {\cal O}_K$. The quotient
$K^*/{\cal O}^*$ will be denoted $\epsilon$. It is a free Abelian
group of rank 2 fitting into a natural exact sequence
$$0\to {\bf Z}\to \epsilon\to {\bf Z}\to 0. \leqno (1.2.1)$$
Let $G$ be as in (1.1). The semidirect
product $\widehat{\widehat W} = 
W\propto (L\otimes\epsilon)$ will be called the double affine Weyl group
for $G$, see [Pa2]. 
Consider the group $G(K)$ and the following three subgroups:
$$D_0 = \{ g\in G({\cal O}'): \pi_k(\pi_K(g))
\in B(F_q)\}, \quad D_1 = \{ g\in G({\cal O}_K): \pi_K(g)\in
T({\cal O}_k) N(k)\}, \leqno (1.2.2)$$
$$D_2 = T({\cal O}') N(K).$$
Thus $D_2$ is a natural ``connected component"
of the Borel subgroup in $G(K)$,
$D_1$ is a similar connected component of the Iwahori subgroup
in $G(K)$ (where $K$ is considered just as a local field with
residue field $k$), and $D_0$ is the "double-Iwahori" subgroup
(cf. [Pa2]).

Let $\varpi_1\in {\bf m}_k$ be a uniformizer of $k$ and $\varpi_2=t$
be the standard uniformizer of $K$. This choice 
of uniformizers gives an identification
$${\bf Z}^2\to\epsilon, \quad (i,j) \mapsto \varpi_1^i \varpi_2^j \,\,
{\rm mod}\,\,({\cal O}')^*$$
and therefore an identification $L\oplus L\to L\otimes\epsilon$,
as well as a realization of $L\otimes\epsilon$ as a subgroup in
$K^*$. Choose also a lifting of $W$ to a subgroup in $G({\cal O}')$.
Then we get an embedding $\widehat{\widehat W}\i  G(K)$.

\proclaim (1.2.3) Proposition.  For any $i,j=0,1,2$ we have the
decomposition
$$G(K) = \coprod_{w\in\widehat{\widehat W}} D_i w D_j$$
and the resulting identification $D_i\backslash\Gamma/D_j\to
\widehat{\widehat W}$ is canonical (independent on the choice of lifting).

\noindent {\sl Proof:} If $E$ is a field, we have the Bruhat decomposition
 $G(E) = \coprod_{w\in W}B(E)wB(E)$. If $(E, \overline E)$
is a local field, then we have the Bruhat-Tits  and Iwasawa decompositions
$$G(E) \quad = \quad \coprod_{w\in \widehat W} IwI \quad =\quad 
\coprod_{a\in L} G({\cal O}_E) a N(E)$$
 where $I$ is the Iwahori
subgroup.
The proposition is obtained by iterated application of the Bruhat decomposition
to the fields $K, k, F_q$ and of the Bruhat-Tits and Iwasawa decompositions
to the local fields $(K,k)$ and $(k, F_q)$, see the argument in [Pa2], \S 2
for the group $PGL_n$. 

\vskip .1cm

In this paper we will be mostly interested in the subgroup $D_1$.

Let $I=\pi_K^{-1}(B(k))$ be the Iwahori subgroup in $G(K)$.
Consider the following homogeneous spaces:
 $$\widehat  X= G(K)/G({\cal O}_K), \quad \widehat  F = G(K)/I_K, \quad
{\cal M}= G(K)/D_1.\leqno (1.2.4)$$
The set
 $\widehat  F$ is the ``affine flag variety" of $G$, see [Lu 1].
Let ${\cal B}(G, K/k)$ be  the Bruhat-Tits building  associated to
$G$ and the local field $(K,k)$. Then $\widehat  X$ is a $G(K)$ -orbit
on the set of
vertices of ${\cal B}(G, K/k)$ while $\widehat  F$ is a $G(K)$-orbit on the
set of flags of type (vertex, maximal cell) in ${\cal B}(G, K/k)$.
According to the general properties of affine buildings, see [Br],
the link of any vertex $v$ of ${\cal B}(G, K/k)$ is a spherical
building associated to $G$ and $k$, and as such, is the boundary
of a locally finite Bruhat-Tits building $\beta_v$ isomorphic
to ${\cal B}(G, k/F_q)$
(conveniently referred to as ``microbuilding"). Cf. the construction of
the double Bruhat-Tits building by A. N. Parshin [Pa2]. 
The set $\cal M$ is naturally identified
with the set of all the horocycles in all the microbuildings $\beta_v$. 
We have the projections
$$p_1: \widehat  F\to\widehat  X, \quad p_2: {\cal M}\to\widehat  F\leqno (1.2.5)$$
with fibers of $p_1$ identified with $F=G(k)/B(k)$ and the fibers of $p_2$
being $L$-torsors.

\vskip .2cm

The Bruhat-Tits decomposition associated to $G$ and the local field
$(K,k)$ allows us to speak about the relative position
$w(b,b')\in\widehat  W$ of two points  $b,b'\in \widehat  F$.
Denote by $U_w(b)\i \widehat  F$ the set of $b'$ such that $w(b,b')=w$
(affine Schubert cell).
 Set also
$$\bar U_w(b) = \bigcup_{w'\leq w} U_w(b), \leqno (1.2.6)$$
where $w'\leq w$ stands for the Bruhat order. For the following,
see [Lu 1].

\proclaim (1.2.7) Theorem. Each $\bar U_w(b)$ has a natural structure
of a projective algebraic variety over $k$ of dimension $l(w)$,
so that $U_w(b)$ is an open subvriety isomorphic to the affine space.
The set $\widehat  F$ has therefore a structure of (the set of $k$-points of)
an ind-object in the category of projective algebraic varieties over $k$.
The group $G(K)$ acts on $\widehat  F$ by automorphisms of an ind-object.

\proclaim (1.2.8) Corollary. $\widehat  F$ has a natural topology, the inductive
limit of compact totally disconnected topologies on the $\bar U_w(b)$.
 The action of $G(K)$ is by homeomorphisms. 

\vskip .2cm

\noindent {\bf (1.3) The central extension.}   Recall the general theory
of Steinberg-Moore-Matsumoto [Mat].
 Let  $E$ be any field, $B$ be an Abelian group and
$s: E^*\times E^*\to B$ be a Steinberg symbol, i.e., a bi-homomorphic
map satisfying the identity $s(x, 1-x)=1$. To this data
Matsumoto associates a central extension
$$1\to B\to \widetilde G(E)_s\buildrel p\over
\to G(E)\to 1.\leqno (1.3.1)$$
As with any central extension, (1.3.1) gives rise to the commutator
map\
$$c_s: G(E)\times G(E)\to B, \quad c_s(x,y) = [\widetilde  x, \widetilde  y]\in B, \,\,
p(\widetilde  x)=x, p(\widetilde  y)=y,\leqno (1.3.2)$$
which is characterized by the property:
$$c_s( \lambda^a, \mu^b) = c(\lambda, \mu)^{ \Psi (a, b)},
 \quad \lambda, \mu\in K^*,
\,\, a,b\in L = {\rm Hom}({\bf G}_m, T).\leqno (1.3.3)$$
Here $\Psi$ is as in (1.1.1).
We specialize this construction to $E=K, B=k^*$ and $s$ being the
tame symbol
$$s_{\rm tame}(x,y) = (-1)^{{\rm ord}(x) {\rm ord}(y)}
\pi\left({x^{{\rm ord}(y)}\over y^{{\rm ord}(x)}}\right).\leqno (1.3.4)$$ 
We denote the corresponding central extension by
$$1\to k^*\to \Gamma \to  G(K)\to 1.\leqno (1.3.5)$$

\proclaim (1.3.6) Proposition. The extension $\Gamma$ canonically
splits over any of the subgroups $D_i\i\Gamma$ defined in (1.2.2).

\noindent {\sl Proof:} Let $M={\rm Norm}(T)$ be the normalizer of $T$.
Matsumoto's explicit construction of $\widetilde G(E)_s$
for any $E,s$, is done in two steps, see [Mat], [Mil].
First, one considers the group scheme $M={\rm Norm}(T)$, the
normalizer of $T$ and constructs a central extension $\widetilde M(E)_s\to M(E)$
by using $s$. Then, one considers the retraction
$\rho: G(E)\to M(E)$ given by the Bruhat decomposition, i.e.,
$\rho(g)$ is the unique element $m\in M(E)$ such that $g=u_1mu_2$
with $u_i\in N(E)$. The extension $\widetilde G(E)_s$ is defined
as a subgroup of bijections of $\widetilde  M(E)\times_{M(E)} G(E)$
generated by certain bijections $u_\alpha(e), e\in E$, $h_\alpha(e),
e\in E^*$, $w_\alpha$,
$\alpha\in R_+$, mimicking  respectively
the action of Chevalley generators from $N(E)$, of the elements of $T(E)$ and
of elements of $W$ labelled by $\alpha$. In our case $E=K, B=k^*$, $s=s_{\rm tame}$. So it is clear that the extension is trivial on $D_2
= N(E)T({\cal O}')$, as $s$ vanishes on
${\cal O}^*\supset ({\cal O}')^*$. As for $D_i, i=0,1$, they both lie in $G({\cal O})$,
and in fact the extension is trivial over $G({\cal O})$ because of
the decomposition $G({\cal O}) = N({\cal O}) M({\cal O})N({\cal O}) $
and again of triviality of the tame symbol on ${\cal O}^*$.

\vskip .1cm

According to the above proposition, we can and will view $D_i$
as subgroups in $\widetilde \Gamma$ and
denote $\Delta_i = {\cal O}_k^*\cdot D_i\i\Gamma$. 
We will be mostly interested in the case $i=1$  and denote
$\Xi =\Gamma/\Delta_1$. Thus
$\Xi\to {\cal M} = \Gamma/\Delta_1$ is a {\bf Z}-torsor, or, equivalently,
$$p: \Xi\to \widehat F\leqno (1.3.7)$$
is an $L_{\rm aff}$-torsor (in the set-theoretical sense). 
For every $b\in\widehat  F$ we denote $\Xi_b = p^{-1}(b)$.

 \vskip .3cm

\noindent {\bf (1.4) The affine Heisenberg-Weyl group.}   
 Note that the composite map
$$K^*\otimes_{\bf Z} K^*\buildrel \{, \,\}\over\rightarrow k^*
\buildrel {\rm ord}_k \over\rightarrow {\bf Z}$$
descends to a nondegenerate skew-symmetric pairing
$\sigma: \Lambda^2(\epsilon)\to {\bf Z}$, where $\epsilon$ is as in (1.2.1).
We get therefore a $W$-invariant
skew-symmetric pairing $\Psi\otimes\sigma$ on $L\otimes\epsilon$ and can form
the corresponding central extension
$$1\to {\bf Z}\to  \widetilde{L\otimes\epsilon} \to L\otimes\epsilon
\to 1,\leqno (1.4.1)$$
(the Heisenberg group) for which the commutator pairing is equal
 to $\Psi\otimes\sigma$. This pairing being $W$-invariant, we have
a natural $W$-action on $\widetilde{L\otimes\epsilon}$. 

\proclaim (1.4.2) Definition. The (double) affine Heisenberg-Weyl group
is the semidirect product $\widetilde W := W\propto (\widetilde {L\otimes
\epsilon})$. 

The following is then straightforward. 

\proclaim (1.4.3) Proposition.  $\widetilde  W$ is isomorphic to $\widehat  W\propto L_{\rm aff}$, the semidirect
product with respect to the action (1.1.3). 

The embedding $\widehat{\widehat{W}}\i G(K)$ chosen in (1.2) induces an embedding
$\widetilde  W\i\Gamma$ and (1.2.3) implies

\proclaim (1.4.5) Proposition.  For any $i,j=0,1,2$ we have the
decomposition
$$\Gamma = \coprod_{w\in\widetilde  W} \Delta_i w \Delta_j$$
and the resulting identification $\Delta_i\backslash\Gamma/\Delta_j\to
\widetilde  W$ is canonical (independent on the choice of liftings).

From this, together with (1.4.3), we deduce:

\proclaim (1.4.6) Proposition. (a) For every $b, b'\in\widehat F$
such that $w(b, b')=w\in\widehat W$, we have a natural identification
of torsors $j_{bb'}: \widehat \Xi_b\to\widehat \Xi_{b'}$, compatible
with the identification $w: L_{\rm aff}\to L_{\rm aff}$ of the structure
groups. These identifications are $\Gamma$-equivariant
and satisfy the transitivity conditions. \hfill\break
(b) The $\Gamma$-orbit on $\Xi\times\Xi$ associated to $(w,l)\in\widetilde
W = \widehat W\propto L_{\rm aff}$, is the subset
$$\Sigma_{w,l} = \biggl\{ (\xi, \xi')\in \widehat\Xi\times\widehat\Xi\biggl|
\xi\in\widehat\Xi_b, \xi'\in\widehat\Xi_{b'}, w(b,b')=w, \xi'=j_{bb'}(\xi)
+l\biggr\}.$$
In particular, $p: \Xi\to\widehat F$ identifies $\Sigma_{w,l}\cap (\{\xi\}
\times\Xi)$
with the Schubert cell $U_w(p(\xi))$.

The purpose of this paper is to make sense of the Hecke algebra
of $\Gamma$ by $\Delta_1$.

\vskip 2cm

\centerline {\bf \S 2. Cherednik algebras.}

\vskip 1cm

\noindent{\bf (2.1) The definitions.} We keep the notation of
(1.1). As $G$ is assumed simply connected, $L$ is spanned by coroots.
Let $Q\i L^\vee$ be the lattice spanned by the roots, so that $L\i Q^\vee$.
Note that $\widehat  W = W\propto L$ while $\widehat  W_{\rm ad} = W\propto Q^\vee$ is 
the extended
affine Weyl group corresponding to $G_{\rm ad}$,
the adjoint quotient of $G$. 
Let $m\in {\bf Z}_+$ be the minimal number such that $m\cdot (a,b)\in
{\bf Z}$ for any $a\in Q^\vee, b\in L^\vee$. Set
 $$P_{\rm aff} = L^\vee\oplus {1\over m} {\bf Z}, \quad \widetilde  T^{\rm aff}
({\bf C}) = {\rm Spec}\, {\bf C}[P_{\rm aff}] = \widetilde  T({\bf C})
\times {\rm Spec} \, {\bf C}[\zeta^{\pm 1/m}].
\leqno (2.1.1)$$
Here $\zeta$ is an independent variable which we will think of
as the second coordinate on the torus $\widetilde  T^{\rm aff}({\bf C})$,
writing a typical point of this torus as $\lambda = (\overline\lambda, \zeta)$
with $\lambda\in\widetilde  T({\bf C})$. The group $\widehat  W_{\rm ad}$ acts on
$\widetilde  T^{\rm aff}({\bf C})$ by
$$a (\lambda,\zeta) = (\zeta^a\cdot\lambda, \zeta), \,\, a\in Q^\vee \quad
w(\lambda, \zeta) = (w(\lambda), \zeta), \,\, w\in W. \leqno (2.1.2)$$
Here the meaning of $\zeta^a$ is as follows. To $a$, there corresponds
a homomorphism $L^\vee\to {1\over m} {\bf Z}$ taking $b\mapsto (b,a)$.
To this homomorphism we associate a homorphism of tori
${\rm Spec} \, {\bf C}[\zeta^{\pm 1/m}] \to \widetilde  T$
whose value on $\zeta$ is denoted by $\zeta^a$.

Let $\Pi \i \widehat  W_{ad}$ be the subgroup of elements of length  0.
It acts naturally on $P_{\rm aff}$. 

Let $r\in {\bf C}^*$ be a fixed nonzero number.The Cherednik algebra
${\cal H}={\cal H}_{r,}$, associated to the root system of $G$ is, by
definition [Ch], generated by the elements
$$\zeta^{\pm 1/m}, \quad \tau_w, w\in \widehat  W^{\rm }, \quad \tau_\pi, \pi\in\Pi, \quad Y_b, b\in P$$
subject to the following relations:

\noindent (2.1.3) $\zeta$ is central and 
the $Y_b$ form the group algebra ${\bf C}[L^\vee]$, i.e.,
$Y_0=1$ and $Y_bY_{b'}=Y_{bb'}$; the $\tau_w$ and $\tau_\pi$ form an affine
Hecke algebra  of $G_{ad}$, see, e.g., [GKV]. We abbreviate $\tau_{s_\alpha}$, 
$\alpha\in R_+$, to $\tau_\alpha$ and $\tau_{s_{\alpha_0}}$ to $\tau_0$.
We also write $Y_{(b,n)} = Y_b \zeta^n$ for $(b,n)\in P_{\rm aff}$.  
$$\tau_\alpha Y_b \tau_{\alpha}^{-1} = Y_b Y_\alpha^{-1}, \quad \alpha\in R_{\rm sim}, (b, \alpha^\vee)=1.\leqno (2.1.4)$$
$$\tau_0Y_b \tau_0^{-1} = Y_b Y_\theta^{-1}\zeta, \quad (b, \theta^\vee)=1.
\leqno (2.1.5)$$
$$\tau_\alpha Y_b = Y_b \tau_\alpha, \quad \alpha\in R_{\rm sim}, 
(b, \alpha)=0.
\leqno (2.1.6)$$
 $$\tau_\pi Y_b \tau_\pi^{-1} = Y_{\pi(b)}, \quad \pi\in\Pi.\leqno (2.1.7)$$

\vskip .2cm

\noindent {\bf (2.2) The residue construction.} In this paper we will
use another construction of the Cherednik algebra which does not use
generators and relations. This construction was given in [GKV]. 
 
 For $\alpha = (n, \bar\alpha)\in\widehat R$, so that $n\in {\bf Z}$ and
$\bar\alpha\in R$, and for $\lambda = (\bar\lambda, n)\in\widetilde  T^{\rm aff}({\bf C})$ we denote
$\lambda^\alpha = \zeta^n \bar\lambda^{\bar\alpha}$ and set, for any
$z\in {\bf C}^*$,
$$\widetilde  T^{\rm aff}_{\alpha, z} = \{\lambda\in\widetilde  T^{\rm aff}({\bf C}):
 \lambda^\alpha = z\}.\leqno (2.2.1)$$
Let ${\bf C}(\widetilde  T^{\rm aff})$ be the field of rational functions on
$\widetilde  T^{\rm aff}({\bf C})$. It is acted upon by $\widehat  W_{ad}$, via
(2.1.2). Let ${\bf C}(\widetilde  T^{\rm aff})[\widehat  W_{ad}]$ be the corresponding
twisted group algebra. Its elements are finite formal sums
$\sum_{w\in\widehat  W_{\rm ad}} f_w(t) [w]$, with $f_w(t)\in {\bf C}
(\widetilde  T^{\rm aff})$
and  the multiplication given by the rules: 
$$[w] [w'] = [ww'], \quad [w] f = f^w [w], \,\, f^w(\lambda) = f(w^{-1}(\lambda)).
\leqno (2.2.2)$$

\proclaim (2.2.3) Theorem. The Cherednik algebra ${\cal H}_{r}$ is
isomorphic to the subalgebra in ${\bf C}(\widetilde  T^{\rm aff})[\widehat  W_{ad}]$ consisting
of $\sum f_w(\lambda)[w]$ such that:\hfill\break
(1) The only possible singularities of each $f_w$ are first order poles
along the $\widetilde  T_{\alpha, 1}^{\rm aff}$, $\alpha\in\widehat  R_+$, with
$${\rm Res}_{\check T^{\rm aff}_{\alpha, 1}}(f_w) +
 {\rm Res}_{\check T^{\rm aff}_{\alpha, 1}}
(f_{s_\alpha w}) = 0, \quad w\in \widehat  W, \alpha\in\widehat  R_+.$$
(2) Each $f_w$ vanishes along each $\widetilde  T^{\rm aff}_{\alpha, r^{2}}$, $\alpha\in 
\widetilde  D(w)$.

\noindent {\sl Proof:}
The statement is almost identical with Theorem 6.3.1 of [GKV] with one
small difference.  Namely, the algebra denoted by $\cher_r$
in [GKV] corresponds to the subalgebra in Cherednik's ${\cal H}_{r}$
generated by $\tau_w, \tau_\pi$ and $Y_b$, where $Y_b$ runs only
over $Q\i L^\vee$. In other words, if we denote by $\check G$ the
Langlands dual group of $G$ and $\check \Pi$ the group of length 0
elements in the extended affine Weyl group of its adjoint quotient, then
${\cal H}_{r} = \cher_r[\check \Pi]$
(the twisted group algebra). Furthermore, the theorem just cited gives
  the description of
$\cher_r$ inside ${\bf C}(\widetilde  T^{\rm aff})[\widehat  W^{}] \i {\bf C}(\widetilde  T
^{\rm aff})
[\widehat  W_{ad}]$ by the same conditions (1) and (2) as in (2.2.3).
It remains to notice that
$${\bf C}(\widetilde  T)[\widehat  W_{ad}] \simeq {\bf C}(\widetilde  T)[\widehat  W^{\rm }][\check \Pi]$$
and that the validity of the conditions (1-2) for an element $\phi$
in the LHS is equivalent to their validity for any coefficient
$\phi_\omega\in {\bf C}(\widetilde  T)[\widehat  W^{\rm }]$ in the decomposition
$\phi=\sum_{\omega\in\check\Pi} \phi_\omega [\omega]$ in the RHS.

\vskip .3cm

\noindent {\bf (2.3) The modified Cherednik algebra.}   
Consider the embedding of lattices
$$L_{\rm aff}  = L\oplus {\bf Z} \buildrel \Psi\oplus
 {\rm Id}\over\hookrightarrow  P_{\rm aff}
=L^\vee\oplus {1\over m} {\bf Z}\leqno (2.3.1)$$
where $\Psi$ is the form (1.1.8),
and the corresponding homomorphism (finite covering) of tori
$$\widetilde  T^{\rm aff} ({\bf C}) = {\rm Spec}\,\, {\bf C}[P_{\rm aff}]
 \to {\rm Spec}\,\,
{\bf C}[L_{\rm aff}] = \check T^{\rm aff}. \leqno (2.3.2)$$
The above maps are $W$-equivariant. Therefore the action (2.1.2) of
$\widehat  W_{ad}$ on $\widetilde  T^{\rm aff}({\bf C})$, being composed of the standard
$W$-action and of the action of $Q^\vee$ by translations, descends to an
action on $\check T^{\rm aff}$. 
Thus  the subfield ${\bf C}(\check T^{\rm aff})\i {\bf C}
(\widetilde  T^{\rm aff})$ is preserved under the
 $\widehat  W_{ad}$-action and can be used
to form a twisted group algebra.

\proclaim (2.3.3) Definition. The modified Cherednik algebra $\cher_{r}(G)$
is the intersection
$${\bf C}(\check T^{\rm aff})[\widehat  W] \cap {\cal H}_{r} \i {\bf C}(\widetilde  T^{\rm aff})[\widehat 
W_{ad}].$$

\vfill\eject

\centerline {\bf \S 3. Hecke operators. Main theorem.}

\vskip 1cm

\noindent {\bf (3.1) Generalities.} 
Let $\Gamma$ be a group,
$\Delta\i\Gamma$ a subgroup. We will denote by $\Gamma//\Delta$
the set of double cosets of $\Gamma$ by $\Delta$. For any such coset
$C=\Delta\gamma\Delta$ we have the Hecke correspondence $\Sigma_C$ which is
the subset
$$\Sigma_C = \bigl\{ (\gamma_1\Delta, \gamma_2\Delta)\bigl| \,\, \gamma_2\gamma_1^{-1}\in C\bigr\} \quad \i\quad (\Gamma/\Delta)\times
(\Gamma/\Delta).\leqno (3.1.1)$$
The $\Sigma_C$ are nothing else than all the $\Gamma$-orbits on
$(\Gamma/\Delta)\times
(\Gamma/\Delta)$.  For $x\in \Gamma/\Delta$ we denote $\Sigma_C(x)$ the
intersection $\Sigma_C\cap ((\Gamma/\Delta)\times \{x\})$.
Denote also 
$$\Gamma/\Delta \buildrel p_1\over\leftarrow \Sigma_C
\buildrel p_2\over\to \Gamma/\Delta\leqno (3.1.2)$$
the natural projections. Thus $\Sigma_C(x) = p_2^{-1}(x)$. 

If $\Gamma$ is a locally compact Hausdorff topological group and $\Delta$
is a compact subgroup, then a Haar measure on $\Gamma$ induces natural 
measures $\mu_{C,x}$ on each $\Sigma_C(x)$ invariant under ${\rm Stab}(x)\i
\Gamma$. Denoting by ${\cal F}(\Gamma/\Delta)$ the space of continuous functions
$\Gamma/\Delta\to {\bf C}$, we have the $\Gamma$-invariant operator
$$\tau_C: {\cal F}(\Gamma/\Delta) \to {\cal F}(\Gamma/\Delta),
\quad (\tau_Cf)(x) = \int_{y\in \Sigma_C(x)} f(y) d\mu_{C,x} \leqno (3.1.3)$$
known as the Hecke operator. Under our assumptions, each $\Sigma_C(x)$
is compact, so the integral makes sense. The Hecke algebra
$H(\Gamma, \Delta)$ is the convolution algebra of continuous compactly
supported, $\Delta$-biinvariant functions on $\Gamma$. Its elements can
be thought of as (continuous) linear combinations of the
operators $\tau_C$. 

\vskip .1cm

We now want to apply the above formalism to the case when $\Gamma$
is the central extension of $G(K)$ introduced in (1.3) and $\Delta =
\Delta_1$ is the ``intermediate Iwahori subgroup". In this case, the
natural topology on $K$ is not locally compact and, moreover, fails to make
it into a topological ring [FP] [Kat], so there is no obvious good topology
on $G(K)$ and $\Gamma$ and no Haar measure. Nevertheless, $\Xi=\Gamma/\Delta_1$
is an $L_{\rm aff}$-torsor over the affine flag variety $\widehat  F=
\bigcup_{w\in\widehat  W} \bar U_w(b_0)$
which has the topology of inductive limit of compact spaces $\bar
U_w(b_0)$. From now on we will freely use the formalism and notation
set up in the Appendix, as well as those of \S 1.

\vskip .3cm

\noindent {\bf (3.2) Function spaces associated to a p-adic loop group.}
The space $\widehat  F$ is an object of $\cal K$ and will also be
regarded as an ind-pro-object of ${\cal S}_0$ (with an action of $G(K)$ by
isomorphisms). 

\proclaim (3.2.1) Proposition. The set-theoretic $L_{\rm aff}$-torsor
$\Xi\to\widehat  F$ has a natural structure of a $\Gamma$-equivariant object
of the category $L_{\rm aff}{\rm -Tors}(\widehat  F)$.

\noindent {\sl Proof:} It is enough to show that for any homomorphism
$\chi: L_{\rm aff}\to {\bf Z}$ the 
induced set-theoretical torsor $\chi_*\Xi\to\widehat  F$ has a natural structure
of a $\Gamma$-equivariant object of {\bf Z}-Tors$(\widehat  F)$. 

Let $\Pi$ be the category of projective algebraic varieties over $k$. Taking
the space of $k$-points defines a functor $\gamma: \Pi\to {\cal P}$. 
The induced functor on ind-objects will be also denoted
$\gamma: {\rm Ind}_s^{\aleph_0}(\Pi) \to {\cal K}$. For a variety $M\in\Pi$
let ${\rm Bun}_1^{\rm alg}(M)$ be the category of algebraic line bundles
on $M$ and for an object $N =(N_i)$ of ${\rm Ind}_s^{\aleph_0}(\Pi)$
we denote ${\rm Bun}_1^{\rm alg}(N) = 2\pro \,{\rm Bun}_1^{\rm alg}(N_i)$. 
For every $M\in\Pi$ we have a natural functor
$${\rm Bun}_1^{\rm alg}(M)\to k{\rm -Bun}_1(\gamma M)$$
(the topological bundle on $k$-points induced by an algebraic bundle).
For every $Y\in {\cal P}$ we have a natural functor
$${\rm ord}_*: k{\rm -Bun}_1(Y) \to {\bf Z}{\rm -Tors}(Y), \quad E\mapsto
(E-\{0\})/{\cal O}_k^*.$$
Composing these functors and passing to ind-objects, we get, for $N\in 
{\rm Ind}_s^{\aleph_0}(\Pi)$, a natural functor
$$\delta: {\rm Bun}_1^{\rm alg}(N) \to {\bf Z}{\rm -Tors}(\gamma N).$$
Now, $\widehat  F$ comes from an object $\underline{\widehat  F}$ of 
${\rm Ind}_s^{\aleph_0}(\Pi)$, in the sense that
 $\widehat  F = \gamma \underline{\widehat  F}$.
Further, it is well known that $\chi$ gives rise to an object
${\cal O}(\chi)\in {\rm Bun}_1^{\rm alg}(\underline{\widehat  F})$,
equivariant under $\Gamma$. It remains to notice that
$\chi_*\Xi = \delta({\cal O}(\chi))$. Proposition is proved.

\vskip .2cm

Thus we can speak about the function spaces ${\cal F}_0(\Xi), {\cal F}(\Xi)$
etc. Let $L^+_{\rm aff}\i L_{\rm aff}$ be the convex cone spanned by the
positive affine roots. For any element $w\in\widehat  W$ let $L^w_{\rm aff} = 
w(L^+_{\rm aff})$. We define ${\cal F}_w(\Xi) = {\cal F}_{L^w_{\rm rat}}(\Xi)$
and similarly define ${\cal F}_w^{\rm rat}(\Xi)$.

\vskip .3cm

\noindent {\bf (3.3) Hecke operators and the main theorems.}
Let $(w,l)\in\widetilde  W = \Gamma//\Delta_1$. see (1.4,6), and $\Sigma_{w,l}
\i\Xi\times\Xi$ be the Hecke correspondence ($\Gamma$-orbit) associated to
$(w,l)$. For $\xi\in\Xi$ let $\Sigma_{w,l}(\xi)\i\Xi$
consist of $\xi'$ such that $(\xi, \xi')\in\Sigma_{w,l}$. 
By (1.4.6), this is a affine space over $k$ of dimension $l(w)$.
The stabilizer ${\rm Stab}(\xi)\i\Gamma$ acts transitively on
$\Sigma_{w,l}(\xi)$.

\proclaim (3.3.1) Proposition. For every $\xi$ the space of complex
valued
Borel measures on $\Sigma_{w,l}(\xi)$ invariant under
${\rm Stab}(\xi)$, is non-zero (and hence is 1-dimensional).

\noindent {\sl Proof:} It is enough to take $\xi=\xi_0$, 
the distinguished point of $\Xi=\Gamma/\Delta_1$, so that
${\rm Stab}(\Xi) = \Delta_1$. In this case the action of
$\Delta_1$ on $\Sigma_{w,l}(\Xi)$ factors through $\Delta_1/\Gamma_0(N)$,
where $N\gg 0$ is a suddificently large integer and
$$\Gamma_0(N) = \{ g\in G({\cal O}_K): g\equiv 1 \,\,\, {\rm mod}\,\,
{\bf m}_K^N\}$$
is the principal congruence subgroup of level $N$. The quotient
$\Delta_1/\Gamma_0(N)$ is an extension
$$1\to E\to \Delta_1/\Gamma_0(N) \to T({\cal O}_k)\to 1,$$
where $E$ is the group of $k$-points of an unipotent algebraic group
over $k$. Further, $\Delta_1/\Gamma_0(N)$ acts on the affine space
$\Sigma_{w,l}(\xi)$ by polynomial transformations, and the action of
$E$ comes from an algebraic action of an algebraic group. It follows
from this and from the compactness of $T({\cal O}_k)$ that for every
$g\in \Delta_1/\Gamma_0(N)$ and any $g$-fixed point $\eta\in\Sigma_{w,l}(\xi)$
the jacobian determinant  $\det(d_\eta g)\in k^*$ has $k$-absolute value 1.
This implies the proposition.

\vskip .2cm

We conlude (by using translation)
that a choice of a ${\rm Stab}(\xi)$-invariant measure
$\mu$ on  $\Sigma_{w,l}(\xi)$ for some one $\xi$ defines unambiguously
 a measure
$\mu_{\xi'}$ on $\Sigma_{w,l}(\xi')$ for any $\xi'$. Notice also that for
any $l, l'\in L_{\rm aff}$ the spaces $\Sigma_{w,l}(\xi)$ and 
$\Sigma_{w,l'}(\xi)$ are canonically identified (1.4.6). 

So for every $w\in\widehat  W$ we choose in some way an invariant measure
$\mu_{w,0,\xi_0}$ on $\Sigma_{w,0}(\xi_0)$ and then define
the measure $\mu_{w,l,\xi}$ on $\Sigma_{w,l}(\xi)$ 
by using the above identifications.
Thus, for a continuous function $f: \Xi\to {\bf C}$ we can formally
write the integral
$$(\tau_{w,l}f)(\xi) = \int_{\eta\in\Sigma_{w,l}(\xi)}
f(\eta) d\mu_{w,l,\xi}\leqno (3.3.2)$$
defining the Hecke operator. Since the domain of integration is
noncompact, the integral may not converge. If $w=e$,
the integration is over a point and $\tau_{e,l}$ is the (well-defined)
operator of shift by $l$. We now state three theorems describing the regularization
of the $\tau_{w,l}$, of their compositions and the structure of the algebra
formed by the regularized operators. The proofs will be given in
the next section. 

\proclaim (3.3.3) Theorem. If $f\in |{\cal F}_0(\Xi)|$, then
${\rm Supp}(f)\cap \Sigma_{w,l}(\xi)$ is compact for any $w,l,\xi$,
so the integral (3.3.2) converges and gives rise to a well-defined
morphism (Hecke operator)
$$\tau_{w,l}: {\cal F}_0(\Xi) \to {\cal F}(\Xi) \quad\in\quad {\rm Mor}(
{\rm Pro}({\rm Mod}_{{\bf C}[L]})).$$
In particular, $\tau_{e,l}$ is the shift by $l$ (preserving 
${\cal F}_0(\Xi)$) and $\tau_{w, l+l'} = \tau_{w,l'}\tau_{e,l}$.

For $l\in L_{\rm aff}$ the corresponding element of ${\bf C}[L_{\rm aff}]$
will be denoted $t^l$, so that a generic element will be written
as a Laurent polynomial $f(t)=\sum_l a_lt^l$. 
The operator $\tau_{w,0}$ will be appreviated to $\tau_w$.
Let $H_{\rm pol}(\Gamma, \Delta_1)$ be the space of formal finite
{\bf C}-linear combinations
$$\sum_{w,l} a_{w,l} \tau_{w,l} = 
\sum_{w\in \widehat  W} f_w(t) \tau_w,\quad a_{w,l}\in {\bf C}, f_w\in 
{\bf C}[L_{\rm aff}]$$
of the Hecke operators corresponding to elements of $\Gamma//\Delta_1$.
This space is not yet an algebra since the $\tau_{w,l}$ act from one
space to another. It is, however, a  ${\bf C}[L_{\rm aff}]$-module
and we get a  ${\bf C}[L_{\rm aff}]$-linear map
$$H_{\rm pol}(\Gamma,\Delta_1) \to {\rm Hom}_{{\rm Pro}({\rm Ind}(
{\rm Mod}_{{\bf C}[L_{\rm aff}]}))}({\cal F}_0(\Xi), {\cal F}(\Xi)).
\leqno (3.3.4)$$
Recall that the torus ${\rm Spec}\, {\bf C}[L_{\rm aff}]$ is denoted
$\check T_{\rm aff}$, so ${\bf C}(\check T_{\rm aff})$ is the field of fractions
of ${\bf C}[L_{\rm aff}]$. 

\proclaim (3.3.5) Theorem. The operator $\tau_{w,l}$ from (3.3.3)
takes values in ${\cal F}_w^{\rm rat}(\Xi)\i {\cal F}(\Xi)$ and gives rise to
an operator 
$$\tau_{w,l}^{\rm rat}\in {\rm End}_{{\rm Pro}({\rm Mod}_{{\bf C}(\check T_{\rm aff})})}({\cal F}^{\rm rat}(\Xi))$$
fitting into the commutative diagram in ${\rm Pro}({\rm Ind}({\rm Mod}_{
{\bf C}[L_{\rm aff}]}))$:
$$\matrix{
&{\cal F}_0(\Xi) & \buildrel \tau_{w,l}
\over\longrightarrow & {\cal F}_w^{\rm rat}
(\Xi)&\cr
c&\bd&&\bd&\Sigma_w\cr
&{\cal F}^{\rm rat}(\Xi)&\buildrel \tau_{w,l}^{\rm rat}\over\longrightarrow
&{\cal F}^{\rm rat}(\Xi)&}$$
where $c$ is the canonical embedding and $\Sigma_w$ is the summation map
from (A.5.2). 

The operator $\tau_w^{\rm rat}$  can be thought of as a regularized Hecke
operator. It now acts from a vector space to itself. This is achieved by
a regularization procedure consisting in summation of a series to a 
rational function and re-expansion in a different domain. So it now makes
sense to ask about the composition of the $\tau_{w,l}^{\rm rat}$.
As before, set $\tau_w^{\rm rat}= {\tau}^{\rm rat}_{w,0}$. Let
$$H_{\rm rat} (\Gamma, \Delta_1) = {\bf C}(\check T_{\rm aff}) \otimes
_{{\bf C}[L_{\rm aff}]} H_{\rm pol}(\Gamma, \Delta_1) = 
\biggl\{ \sum_{w\in\widehat  W} f_w(t) \tau_w^{\rm rat}\biggl|
f_w\in {\bf C}(\check T_{\rm aff})\biggr\}
\leqno (3.3.6)$$
be the space of formal ${\bf C}(\check T_{\rm aff})$-linear combinations 
of the $\tau_w^{\rm rat}$. The map (3.3.4) induces a 
${\bf C}(\check T_{\rm aff})$-linear map
$$\tau: H_{\rm rat}(\Gamma, \Delta_1) \to {\rm End}_{{\rm Pro}({\rm Mod}_
{{\bf C}[L_{\rm aff}]})} ({\cal F}^{\rm rat}(\Xi)).\leqno (3.3.7)$$
Expanding the rational functions $f_w$ into power series in some domaiin,
we can view elements of  $H_{\rm rat}(\Gamma, \Delta_1)$ as certain 
{\it infinite} formal linear combinations of the operators corresponding to 
double cosets. It turns out that considering such combinations is necessary for
the algebra generated by the $\tau_w^{\rm rat}$ to close. 

\proclaim (3.3.8) Theorem. (a) The space $H_{\rm rat}(\Gamma, \Delta_1)$
has a natural structure of an associative algebra so that
$\tau$ is a homomorphism of algebras.\hfill\break
(b) Let $H(\Gamma, \Delta_1)\i H_{\rm rat}(\Gamma, \Delta_1)$
be the subalgebra formed by $S$ such that $\tau(A)$ preserves
${\cal F}_0(\Xi) \i {\cal F}^{\rm rat}(\Xi)$. Then  $H(\Gamma, \Delta_1)$
is isomorphic to the modified Cherednik algebra $\cher_q(G)$ from
(2.3.3).

\vskip .1cm

\noindent {\bf (3.3.9) Remarks.} (a) The above approach can be
applied
to the case of the locally compact group $G(k)$ and the non-compact
subgroup $N(k) T({\cal O}_k)$. In this case the properties of
p-adic intertwiners together with Theorem 5.4 of [GKV] can be
used to give  a  simple proof of the main result of [CK].

(b) The shape of the particular infinite combinations of double cosets
appearing in the Hecke algebra $H(\Gamma, \Delta_1)$ and
in the simpler situation of Remark (a),  is remindful of
the procedure of perverse (or  intersection homology) extension
of sheaves. The right context here seems to be the
(not yet developed) theory of ``semi-infinite'' perverse
extension as sketched by Lusztig [Lu 2].

\vskip .1cm

\proclaim (3.3.10) Corollary. For any subgroup $U\i \Gamma$ the
(pro-)space of invariants ${\cal F}_0(\Xi)^U$ is naturally acted
upon by $\cher_q(G)$. 

\vskip .1cm

\noindent{\bf (3.3.11) Example.} An interesting example of a subgroup
$U\i\Gamma$ is obtained as follows. Let $C$ be a smooth irreducible
projective curve over $k$ and $x\in C$ be a $k$-point. Let $K=k(C)_x$
be the completion of $k(C)$ at $x$. Then $K\simeq k((t))$, so
(3.3.9) is applicable to $U:= G(k[C-\{x\}])\i\Gamma$. We get
an action of $\cher_q(G)$ on $|{\cal F}_0(\Xi)|^U$. Elements of
the latter space can be viewed as certain functions on isomorphism
classes of principal $G$-bundles on $C$ equipped with a 
``horocyclic structure'' (i.e., with a reduction of the structure
group from $G(k)$ to $N(k) T({\cal O}_k)$) at $x$ and with 
a choice of a nonzero vector in the determinantal space
$\det \, H^\bullet(C, {\rm ad}(P))$. We get in this way a setup
for generalizing the theory of automorphic forms over function
fields [Dr] from the case of a finite to the case of a p-adic field
of constants.  Hecke operators are now (linear combinations of)
singular integral operators, while in the classical case they are
given by finite summation.

\vskip 2cm

\centerline {\bf \S 4. Principal series representations and the proof
of main theorems.}

\vskip 1cm

\noindent {\bf (4.1) The unramified principal series.} We  specialize the
 discussion of (A.7)
 to 
$$A=L_{\rm aff}, \quad T_A = \check T_{\rm aff}, 
\quad Z=\widehat  F, \quad B=\Xi,$$
see (3.2).
 Because $\Xi$ is a $\Gamma$-equivariant $L_{\rm aff}$-torsor
over $\widehat  F$, we get $\Gamma$-equivariant line bundles ${\cal L}(\lambda)$,
$\lambda\in\check T_{\rm aff}$, on $\widehat  F$. The ``space" of sections
$V_\lambda = \Gamma(\widehat  F, {\cal L}(\lambda))$ is acted upon
by $\Gamma$. It will be called the {\it unramified principal
series}  representation of $\Gamma$. According to our general principles,
we consider it as an object of ${\rm Pro}({\rm Vect})$.

\vskip .1cm

\noindent {\bf (4.1.0) Remark.} At this point it is natural to ask
what is the analog, for groups such as $\Gamma$, of the concept
of a smooth representation of a p-adic group
(so that our $V_\lambda$ are examples). 
 While this will not be used in this paper, let us indicate the answer. 
The field
$K$ is naturally (the set of maps from $\{ {\rm pt}\}$ to)
a ring object in the category ${\rm Pro}({\rm Ind}({\rm Pro}({\rm Ind}
({\cal S}_0))))$, see [Kat]. The group $\Gamma$ is naturally a group
object in this category. On the other hand, formulas for
Hom-sets in ind- and pro-categories show that for a
pro-{\bf C}-vector space $E$ (i.e,. a pro-ind-object in ${\rm Vect}_0$)
the endomorphism algebra ${\rm End}(E)$ is naturally a ring object of the
two-fold iterated pro-ind-category of affine algebraic varieties over
{\bf C}, and  ${\rm Aut}(E)$ is a group object
in this category. 
Replacing algebraic varieties by their sets of points, let us
consider ${\rm Aut}(E)$ as a group 2-fold iterated pro-ind-object
in $\cal S$. A (2-)smooth representation is, by definition,
 just a morphism of such group objects.

\vskip .2cm

 The pro-${\bf C}(\check T
_{\rm aff})$-vector space $\Gamma_{\rm rat}(\check T_{\rm aff}, V)$ will
be called the {\it generic principal series} representation. The Mellin
transform from (A.7) gives identifications
$${\cal F}_0(\Xi)\to\Gamma_{\rm reg}(\check T_{\rm aff}, V), \quad 
{\cal F}^{\rm rat}(\Xi) \to \Gamma_{\rm rat}(\check T_{\rm aff}, V). \leqno
(4.1.1)$$
For $b\in \widehat  F$ let $I_b\i \Gamma$ be the stabilizer of $b$ and $N_b=[I_b, I_b]$. For $w\in\widehat  W$ denote by $\mu_w(b)$ the vector space of
$N_b$-invariant {\bf C}-valued Borel measures on the Schubert cell
$U_w(b)$. By the same reason as in (3.6.1), we have $\dim_{\bf C}(\mu_w(b))
= 1$.

\proclaim (4.1.2) Proposition. (a) When $w\in\widehat  W$ is fixed, the
$\mu_w(b)$ form the fibers of a $\Gamma$-equivariant line bundle
$\mu_w$ on $\widehat  F$ in the sense of (A.4.2).\hfill\break
(b) The bundle $\mu_w$ is $\Gamma$-equivariantly isomorphic to
${\cal L}(q^{\widehat \delta_w})$, where $\widehat \delta_w$ was defined in
(1.1.6).

\noindent {\sl Proof:} (a) Let $U_w\in\widehat  F\times\widehat  F$ be the
Schubert correspondence and $p_w: U_w\to \widehat  F$ be the projection
onto the second factor. The fact the the $\mu_w(b)$ tie together to
form a line bundle on $\widehat  F$ in the sense of ind-objects in
${\rm Pro}({\cal S}_0)$, follows from the fact that $p_w$ gives rise
to a locally trivial algebraic affine bundle over $\widehat  F$ in the sense
of ind-objects in the category of projective $k$-varieties. This latter
fact is well known. 

(b) It is enough to consider the distinguished point $b_0\in \widehat  F$ with
stabilizer $I_K$. Then, $U_w(b_0)$ has a unique $I_K$-fixed point 
$b_w$ and the eigenvectors of the $I_K$-action on $T_{b_w}U_w(b_0)$ are
in bijection with positive affine roots from $D(w)$, whence the
statement. 

\vskip .2cm

We now proceed to define the analogs, in our affine situation, of the
principal series intertwiners for p-adic groups  [GGP] [Cas].
Consider the action (1.1.3) of $\widehat W$ on $L_{\rm aff}$;
the induced action on $\check T_{\rm aff} = {\rm Spec} \, {\bf C}[L_{\rm aff}]$
will be denoted by $\lambda\mapsto w(\lambda)$. 
Introduce the twisted
$\widehat  W$-action on $\check T_{\rm aff}$ by
$$w*\lambda = q^{\widehat  \rho} \cdot w(q^{-\widehat  \rho}\cdot\lambda) = q^{\widehat  \delta_w}\cdot
w(\lambda) = w(q^{-\widehat  \delta_w}\cdot\lambda). \leqno (4.1.3)$$
Define the linear operator
$$M_w(\lambda): V_\lambda = \Gamma (\widehat  F, {\cal L}_{\lambda q^{-\widehat \delta_w}}
\otimes\mu_w) \to \Gamma (\widehat  F, {\cal L}_{w*\lambda}) = V_{w*\lambda}
\leqno (4.1.4)$$
by
$$M_w(\lambda)(f\otimes m)(\xi) = \int_{b'\in U_w(b)} f(j_{bb'}(\xi))dm(b'), \quad\xi\in\Xi_b \leqno (4.1.5)$$
(provided the integral converges).
Note that up to the choice of the class of functions considered, (4.1.5)
is identical with the formula (3.3.2) defining the Hecke operator $\tau_{w,0}$.

\proclaim (4.1.6) Proposition. 
 (a) If $|\lambda^{\alpha^\vee}|>q^{-1}$ for any
$\alpha\in D(w)$, then the integral (4.2.4) converges and gives rise to
a $\Gamma$-equivariant
morphism $M_w(\lambda): V_\lambda\to V_{w*\lambda}$
in ${\rm Pro}({\rm Vect})$. \hfill\break
(b) The operators $M_w(\lambda)$ extend to a rational (A.6.2) isomorphism $M_w: V\to
w^*V$ of pro-vector bundles on $\check T_{\rm aff}$.\hfill\break
(c) The operators
$$A_w = {1\over c_w(\lambda)} M_w, \quad c_w(\lambda) = \prod_{\alpha\in
D(w)} {1-\lambda^{\alpha^\vee}\over 1-q\lambda^{\alpha^\vee}},$$
satisfy the conditions $A_wA_{w'}=A_{ww'}$, in particular, they form a
representation of $\widehat  W$ in $V_{\rm rat} = 
\Gamma_{\rm rat}(\check T_{\rm aff}, V)$. \hfill\break
(d) The only singularities (A.6.3) of $A_w$ are first order poles at the hypersurfaces
$\check T_{\alpha, 1}^{\rm aff}$, $\alpha\in D(w)$.

\noindent {\sl Proof:} This is similar to the proof of the analogous
statement (see [GGP] [Cas]) 
for the locally compact group $G(k)$ instead of the affine
group $\Gamma$.

First, we consider the case of a simple affine reflection $w=s_\alpha, 
\alpha\in\widehat  R_{\rm sim}$. Let $P_\alpha$ be the parahoric subgroup
in $\Gamma$ corresponding to the set of roots
$\widehat  R_+ \cup \{-\alpha\}$, so that we have the projection
$$p_\alpha: \widehat  F\to \widehat  F_\alpha := 
\Gamma/P_\alpha \leqno (4.1.7)$$
with fibers isomorphic to $P^1(k)$. The Schubert correspondence
$U_\alpha = U_{s_\alpha}\i \widehat  F\times\widehat  F$ consists in this case of
$(b, b')$ such that $p_\alpha(b)\neq p_\alpha(b')$. We can think of the 
ind-scheme structure on $\widehat  F$ to be given in the form $\widehat  F = 
``\ind"_{i\in I} \,p_\alpha^{-1}(\widehat  F_{\alpha, i})$, where
$\widehat  F_{\alpha} = ``\ind" \widehat  F_{\alpha, i}$ is the exhaustion of $\widehat  F
_\alpha$ by its closed Schubert varieties. So the action of $M_{s_\alpha}$
is fiberwise with respect to $p_\alpha$. the situation on each fiber of
$p_\alpha$ is identical to the situation for the intertwining operator for 
the group $SL_2(k)$, for which $P^1(k)$ is the flag variety. So the
well known properties of $SL_2(k)$ imply all the properties of
$M_{s_\alpha}$ claimed in the proposition, including the equality
$A_{s_\alpha}^2=1$. 

Next, let $w=s_{\alpha_1} ... s_{\alpha_n}$, $n=l(w)$, be a reduced
decomposition of an arbitrary $w\in\widehat  W$. Then the properties of the
BN-pair associated to the affine root system, imply that
$$U_w = U_{s_{\alpha_1}} \times_{\widehat  F} U_{s_{\alpha_2}} \times_{\widehat  F} ...
\times_{\widehat  F} U_{s_{\alpha_n}}$$
and thus $M_w = M_{s_{\alpha_1}} \circ ... \circ M_{s_{\alpha_n}}$ is indeed
a rational isomorphism and $A_w$ satisfies the claimed condition
on singularities. Finally, if $w=w_1w_2$ with $l(w)=l(w_1)+l(w_2)$, then
$U_2=U_{w_1}\times_{\widehat  F} U_{w_2}$, so $M_w = M_{w_1}\circ M_{w_2}$ and
$A_w = A_{w_1}\circ A_{w_2}$. This, together with $A_{s_\alpha}^2=1$,
$\alpha\in\widehat  R_{\rm sim}$, implies that $A_{w_1w_2}
=A_{w_1}\circ A_{w_2}$ for any $w_1, w_2$. Proposition is proved. 

\vskip .3cm

\noindent {\bf (4.2) Proof of Theorem 3.3.3.} We denote $p: \Xi\to\widehat  F$
the projection. Let us equip $\Xi$ with the topology of inductive limit
of locally compact spaces ($L_{\rm aff}$-torsors) 
$\Xi^{(\gamma)} =p^{-1}(\bar U_\gamma(b_0))$, $\gamma\in\widehat  W$.
Then any $\Sigma_{w,l}(\xi)$ is contained in some $\Xi^{(\gamma)}$. 
Every $f\in |{\cal F}_0(\Xi)|$ gives, upon restriction to any
$\Xi^{(\gamma)}$, a function which is compactly supported in the
topological sense.  This the convergence of $(\tau_{w,l}f)(\xi)$
would follow from the next lemma.

\proclaim (4.2.1) Lemma. Each $\Sigma_{w,l}(\xi)$ us closed in $\Xi$
(and hence in any $\Xi^{(\gamma)}$ containing it).

\noindent {\sl Proof:} It is enough to consider the case $l=0$.
Denote $\Sigma_{w,0}$ by $\Sigma_w$. If $w= s_{\alpha_1} ... s_{\alpha_n}$,
$\alpha_i\in\widehat  R_{\rm sim}$, as a reduced decomposition, then
$$\Sigma_w = \Sigma_{s_{\alpha_1}} \times_\Xi \Sigma_{s_{\alpha_2}}
\times_\Xi ... \times_\Xi \Sigma_{s_{\alpha_n}}.$$
Suppose that we know that all the $\Sigma_{s_{\alpha_i}}(\xi)$ are closed.
Then $\Sigma_w(\xi)$ is obtained by taking the closed subset
$\Sigma_{s_{\alpha_1}}(\xi)$, for any point $\xi_1$ of this subset,
taking the closed subset $\Sigma_{s_{\alpha_2}}(\xi_1)$, and so on,
so $\Sigma_w(\xi)$ will be closed.

So we reduce to the case $w=s_\alpha$, $\alpha\in\widehat  R_{\rm sim}$.
In this case the statement reduces to one about the group $SL_2(k)$ and
its homogeneous space
$$\Xi_0 = SL_2(k)/B^0, \quad B^0 = \pmatrix{ {\cal O}_k^*& k\cr 0& {\cal O}_k^*}
\cap SL_2(k).$$
we have to prove that $B^0$-orbits in $\Xi_0$ are closed. But $\Xi_0 = 
(k^2-\{0\})/{\cal O}_k^*$, and $B^0$-orbits not reducing to single points,
are the images of straight lines in $k^2-\{0\}$, so they are closed. 
Lemma is proved.

\vskip .2cm

Having established the lemma, we get the convergence of $\tau_{w,l}(f)$
at the level of functions, and a standard argument (cf. the proof of
(4.1.2)(a)) shows that in this was we get a morphism of pro-objects as claimed
in our theorem. The rest of the claims of the theorem (about $\tau_{e,l}$)
are clear. 

\vskip .3cm

\noindent {\bf (4.3) Proof of Theorem 3.3.5.} We first prove
that $\tau_{w,l}({\cal F}_0(\Xi)) \i {\cal F}_w(\Xi)$. For this, it is
enough to show the following.

\proclaim (4.3.1) Lemma. Let $\xi\in\Xi$, $b=p(\xi)\in\widehat  F$. Let
 $\gamma\in\widehat  W$ be such that $U_w(b)\i \bar U_\gamma(b_0)$. 
Take a continuous section $s$ of the $L_{\rm aff}$-torsor
$p^{-1}(\bar U_\gamma(b_0))\to \bar U_\gamma(b_0)$ 
(i.e., of the restriction of $\Xi$) and write, for any $b'\in U_w(b)$,
$$j_{bb'}(\xi) = s(b') + a(b')$$
so that $a: U_w(b)\to L_{\rm aff}$ is a locally constant function.
Then the image of $a$ is contained in some affine translation of
$w(L^+_{\rm aff})$. 

\noindent {\sl Proof:} As with some of the previous lemmas, this statement
reduces to the case when $w$ is a simple affine reflection, i.e., to
a similar statement about the {\bf Z}-torsor $(k^2-\{0\})/{\cal O}_k^*
\to P^1(k)$, which is elementary and left to the reader. 

\vskip .2cm

To prove the rest of the assertions of the theorem, we identify, via the Mellin
transform (A.7.7)
$${\cal F}_0(\Xi) \simeq\Gamma_{\rm reg}(\check T_{\rm aff}, V), \quad 
{\cal F}^{\rm rat}(\Xi) \simeq \Gamma_{\rm rat}(\check T_{\rm aff}, V).$$
After that, we notice that the Hecke operator
$$\tau_{w,0}: {\cal F}_0(\Xi) \to {\cal F}_w(\Xi)$$
is given by the same formula (3.3.2)=(4.1.5) as the operator induced by
the principal series intertwiner $M_w$ on regular rections:
$$\Gamma_{\rm reg}(M_w): \Gamma_{\rm reg}(\check T_{\rm aff}, V) \to
\Gamma_{\rm rat}(\check T_{\rm aff}, W^*V) = 
\Gamma_{\rm rat}(\check T_{\rm aff}, V).$$
We conclude that  $\tau_{w,0}({\cal F}_0(\Xi))\i {\cal F}_w^{\rm rat}(\Xi)$
and define the operator $\tau_{w,0}^{\rm rat}$ to be the operator induced by
$M_w$ on rational sections. Then, we define $\tau_{w,l}^{\rm rat}
=\tau_{0,l}\circ \tau_{w,0}^{\rm rat}$, where $\tau_{0,l}$ is the
operator of shift by $l$. All the statements of the theorem follow
now from the properties of the Mellin transform and of the intertwiners
$M_w$. 

\vskip .3cm

\noindent {\bf (4.4) Proof of Theorem 3.3.8.} (a) Because
we can identify $\tau_w^{\rm rat}$ with the intertwiner $M_w$, an element
of $H_{\rm rat}(\Gamma, \Delta_1)$ can be viewed as a finite
formal linear combination $\sum_{w\in\widehat  W} f_w(t) M_w$ with
$f_w\in {\bf C}(\check T_{\rm aff})$. Since $A_w$ is just a rational
multiple of $M_w$, we can as well think that $H_{\rm rat}(\Gamma, \Delta_1)$
consists of formal linear combinations $\sum f_w(t) A_w$.
Because the $A_w$ form a representation of $\widehat  W$ in 
${\cal F}^{\rm rat}(\Xi) = \Gamma_{\rm rat}(\check T_{\rm aff}, V)$ and take
$V_\lambda\to V_{w*\lambda}$, a natural algebra structure on
$H_{\rm rat}(\Gamma, \Delta_1)$ is obtained by identifying it with the
twisted group algebra ${\bf C}(\check T_{\rm aff})[\widehat  W]_*$ formed
 with
respect to the twisted action of $W$ on $\check T_{\rm aff}$. In other
words, ${\bf C}(\check T_{\rm aff})[\widehat  W]_*$ consists of formal
finite sums
$$\sum_{w\in\widehat  W} f_w(t) [w], \quad [w]\cdot [w'] = [ww'], \quad
 [w]f=f_*^w[w], \quad f_*^w(t):= f(w^{-1}*t).$$
The identification ${\bf C}(\check T_{\rm aff})[\widehat  W]_*\to 
H_{\rm rat}(\Gamma, \Delta_1)$  is just $[w]\to A_w$. Proposition 4.1.6
impies that in this way the map $\tau$ from (3.3.7) becomes a
homomorphism of algebras. 

\vskip .1cm

(b) In view of Definition 2.3.3, Theorem 2.2.3 and of the difference of the
two actions of $\widehat  W$ on $\check T_{\rm aff}$ (the ``straight" action defining
${\bf C}(\check T_{\rm aff})[\widehat  W]$ and the action (4.1.3) defining
${\bf C}(\check T_{\rm aff})[\widehat  W]_*$), it is enough to prove the
following.

\proclaim (4.4.1) Proposition. In order that $\phi=\sum_{w\in \widehat  W}
f_w(t)A_w\in H_{\rm rat}(\Gamma, \Delta_1)$ 
belong to $H(\Gamma, \Delta_1)$, it is necessary and sufficient
 that the following hold:\hfill\break
(1) The only possible singularities of each $f_w$ are poles of order $\leq 1$ 
on
$\check T_{\alpha, q^{-1}}^{\rm aff} , \alpha\in\widehat   R_+$, with
$${\rm Res}_{\check T^{\rm aff}_{\alpha, q^{-1}}}(f_w) + 
{\rm Res}_{\check T^{\rm aff}_{\alpha, q^{-1}}}(f_{s_\alpha w}) = 0, 
\quad w\in \widehat  W, \alpha\in\widehat   R_+.
$$
(2) Each $f_w$ vanishes along each $\check T^{\rm aff}_{\alpha, 1}$, 
$\alpha\in D(w)$.

\noindent {\sl Proof:} Let us summarize the relevant properties of the
$A_w$ which follow from their defininition and from the properties of
the $M_w$ (Proposition 4.1.6):

\vskip .1cm

\noindent (4.4.2) Each $A_w$ has a first order pole along each
$\check T^{\rm aff}_{\alpha, 1}$, $\alpha\in D(w)$, and no other
singularities.

\vskip .1cm

Note that $s_\alpha *\lambda=\lambda$ for
 $\lambda\in\check T^{\rm aff}_{\alpha, q^{-1}}$, so if $\lambda$ is
a generai point of $ T^{\rm aff}_{\alpha, q^{-1}}$, then
$A_{s_\alpha}: V_\lambda\to V_\lambda$ is a well-defined
involution (in the category of pro-vector spaces). 

\proclaim (4.4.3) Lemma. In the above situation,
 $A_{s_\alpha}: V_\lambda\to V_\lambda$
is the identity. 

\noindent {\sl Proof:} If $\alpha$ is a simple affine root, the statement
reduces to that about the intertwiner $A$ in the unramified principal
series of $SL_2(k)$. Namely, that the action of $A$ in the space of
half-measures on $P^1(k)$ (i.e., of sections of the standard square root
of the bundle of measures) is the identity. This is well known 
as the corresponding representation is irreducible [GGP].
The case of a general $\alpha\in\widehat  R_+$ is obtained from this
by conjugation: we take $w\in\widehat  W$ such that $w(\alpha)=
\beta\in\widehat R_
{\rm sim}$, then $A_w$ intertwines the $A_{s_\alpha}$-action on $V_\lambda$,
$\lambda\in \check T^{\rm aff}_{\alpha, q^{-1}}$ and the $A_{s_\beta}$-action
on $V_\mu$, $\mu=w*\lambda\in \check T^{\rm aff}_{\beta, q^{-1}}$. Lemma
4.4.3 is proved.

\vskip .2cm

We need another lemma to pinpoint the origin of the conditions (1) of
Proposition 4.4.1. Let $E$ be a free ${\bf C}[[x]]$-module (possibly of
 infinite
rank). We denote $\bar E=E/xE$ the vector space which can be thought 
of as the fiber of $E$ (as a bundle) at $x=0$. If $E=(E_i)_{i\in I}$ is
a strict pro-object in the category of free ${\bf C}[[x]]$-modules,
we will denote by $\bar E$ the pro-vector space $(\bar E_i)_{i\in I}$.
Let now ${\bf Z}/2 = \{1,s\}$ act on ${\bf C}[[x]]$ by $s^*(x)=-x$
amd suppose that $E$ is a ${\bf Z}/2$-equivariant (pro-) free
${\bf C}[[x]]$-module.
In other words, we have a {\bf C}-linear involution $s: E\to E$ such that
$s(xe)=-xs(e)$. Then $s$ induces an involution $\bar s:\bar E\to\bar E$
of (pro-) vector spaces. 

\proclaim (4.4.4) Lemma. Suppose, in the above situation, that
$\bar s={\rm Id}$. Consider ${\bf C}((x))[{\bf Z}/2]$, the twisted
group algebra with coefficients in ${\bf C}((x))$. It acts naturally
on $E\otimes {\bf C}((x))$. In order that the action of
$f_0(x) + f_1(x) [s]$ preserve $E$, 
 it is necessary and sufficient that the following hold:\hfill\break
(a) ${\rm ord}_x(f_i) \geq -1, \,\,\, i=0,1$. \hfill\break
(b) ${\rm Res}(f_0) + {\rm Res}(f_1)=0.$ 

The proof is left to the reader. 

\vskip .2cm

We now return to the proof of Proposition 4.4.1. Assume that
$\Sigma_w f_w(t) A_w$ preserves the regular sections of $V$.
 As each $A_w$ does indeed have a pole along each $\check T^{\rm aff}_
{\alpha, 1}$, $\alpha\in D(w)$, the condition (2) follows. Since the
$A_w$ have no other singularities, each $f_w(\lambda$ should be regular at
every $\lambda$ whose stabilizer in $\widehat  W$ is trivial, i.e.,
$\lambda\notin \bigcup_{\alpha\in\widehat  R_+} \check T^{\rm aff}_{\alpha,
q^{-1}}$. Further, by taking a generaic point, say, $\lambda_0$, of
some $\check T^{\rm aff}_{\alpha,
q^{-1}}$, we have ${\rm Stab}(\lambda_0)=\{1, s_\alpha\}$. By taking a
1-dimensional formal transversal slice to $\check T^{\rm aff}_{\alpha,
q^{-1}}$ al $\lambda_0$, we see that the rest of the conditions
required in (1), follow from Lemmas 4.4.4 and 4.4.3. The
proof of sufficiency of (1) and (2) for $\sum f_w(t)A_w$
to preserve regular sections, is similar.

This finishes the proof of Proposition 4.4.2 and Theorem 3.3.8.

\vfill\eject

\centerline {\bf Appendix. Analysis on ind/pro objects.}

\vskip 1cm

\noindent {\bf (A.1) Instances of elementary description.}
For a category $\cal C$ we note by ${\rm Ind}({\cal C})$ and ${\rm Pro}
({\cal C})$ the categories of ind- and pro-objects in $\cal C$.
Their objects are filtering inductive or projective systems
 $(X_i)_{i\in I}$ over $\cal C$, i.e., co- or contravariant functors 
$I\to{\cal C}$,
where $I$ is a small filtering category. For background, see [AGV] Exp. I
and [AM]. 
Sometimes we will use the notations $``\ind" X_i$, $``\pro" X_i$ for
the ind/pro-object $(X_i)$. 
By ${\rm Ind}^{\aleph_0}({\cal C})$, ${\rm Pro}^{\aleph_0}({\cal C})$
we will denote the full subcategories formed by $(X_i)_{i\in I}$ with
${\rm Ob}(I)$ countable. By ${\rm Ind}_s({\cal C}$, ${\rm Pro}_s({\cal C})$
we will denote the full subcategories formed by strict ind- or pro-objects,
i.e., by inductive (projective) systems in which all the structure
maps are monomorphisms (epimorphisms).

We will need notations for the following categories (to be used in the
main text):

\vskip .1cm

\noindent ${\cal S}_0, {\cal S}$: finite sets, all sets. 

\vskip .1cm

\noindent ${\rm Vect}_0$, Vect: finite-dimensional {\bf C}-vector
spaces, all vector spaces.

\vskip .1cm

\noindent ${\rm Mod}^0_R, {\rm Mod}_R$ ($R$ a ring): finitely presented left $R$-modules,
all left $R$-modules.

\vskip .1cm

\noindent ${\rm Coh}(S), {\rm QCoh}(S)$ ($S$ an algebraic variety over 
a field): coherent, quasicoherent sheaves of ${\cal O}_S$-modules. 

\vskip .1cm

\noindent $\cal P$: compact Hausdorff totally disconnected spaces.

\vskip .1cm

\noindent ${\cal K}e$: Kelley spaces [GZ], i.e., Hausdorff spaces
$Z$ such that $U\i Z$ is open iff $U\cap C$ is open in $C$ for any compact
$C\i Z$.

\vskip .1cm

\noindent $\cal K$: Kelley spaces which can be represented as a countable
unions of compact subspaces lying in $\cal P$. 

\vskip .1cm

The following is then elementary.

\proclaim (A.1.1) Proposition. (a) The functors of taking the
inductive or projective limit establish equivalences:
$${\rm Ind}({\cal S}_0)\simeq {\cal S}, {\rm Ind}({\rm Mod}^0_R)\simeq
{\rm Mod}_R, {\rm Pro}({\cal S}_0)\simeq {\cal P}, 
{\rm Ind}_s^{\aleph_0}({\cal P})\simeq {\cal K}.$$
(b) The functor $\pro: {\rm Pro}_s({\rm Vect})\to {\rm Vect}$ is injective
on Hom-sets. 

In general, for an ind- or pro-object $X$ we will denote by $|X|$ the
result of actually taking the limit of the inductive/projective
system $X$. For example, if $X=(X_i)$ is a pro-object in ${\cal S}_0$,
we will denote $|X|=\pro X_i$ the corresponding object of $\cal P$ etc. 

\vskip .3cm

\noindent {\bf (A.2) Function spaces.}
For a finite set $X$ let ${\cal F}(X)$ be the space of functions
$X\to {\bf C}$. This gives a contravariant functor
${\cal F}: {\cal S}_0\to {\rm Vect}_0$. We then extend it to pro-objects
componentwise getting a functor, also denoted $\cal F$:
$${\cal F}: {\rm Pro}({\cal S}_0) \simeq {\cal P}\to {\rm Ind}({\rm Vect}_0)
\simeq {\rm Vect.}\leqno (A.2.1)$$
For $Y\in {\rm Pro}({\cal S}_0)$ the space $|{\cal F}(Y)|$ consists of
locally constant functions $|Y|\to{\bf C}$, i.e., functions which ar continuous
if {\bf C} is given the discrete topology. We then extend $\cal F$
to ind-objects componentwise, getting
$${\cal F}: {\cal K}\simeq {\rm Ind}_s^{\aleph_0}({\cal P}) \to {\rm Pro}
_s^{\aleph_0}({\rm Vect}). \leqno (A.2.2)$$ 
Again, for $Z\in {\rm Ind}_s^{\aleph_0}({\cal P})$ the space
$|{\cal F}(Z)|$ consists of continuous functions $|Z|\to {\bf C}$
(with {\bf C} discrete). But now the pro-object ${\cal F}(Z)$ contains
more information than the vector space $|{\cal F}(Z)|$ (namely, the
natural topology of projective limit of discrete spaces). In the
present paper the function space of a space from $\cal K$ will
be always understood as a pro-object.

\vskip .3cm

\noindent {\bf (A.3) 2-limits.} We recall some general categorical
constructions. 

Let $({ C}_i)_{i\in I}$
be a filtering inductive system of categories, with
structure functors denoted $C(\alpha): C_i\to C_j$ for $\alpha\in {\rm Hom}_I
(i,j)$. The direct 2-limit $2\ind_{i\in I} C_i$ is the category whose objects
are pairs $(i, X_i)$ with $i\in I, X_i\in C_i$, with
$${\rm Hom}((i, X_i), (j, X_j)) = \ind_{(k, \alpha, \beta)\in I/\{i,j\}}
{\rm Hom}_{C_k}(C(\alpha)(X_i), C(\beta)(X_j)).\leqno (A.3.1)$$
Here $I/\{i,j\}$ is the category whose objects are diagrams in $I$ of the
form
$$i\to k\leftarrow j$$
with the obvious concept of morphisms.
This category is a certain localization of the category called
the quasi-colimit in [Gra], p. 201 and the homotopy limit in [Th].

If $(C_i)_{i\in I}$ is a filtering projective system of categories,
with structure functors \hfill\break
 $C(\alpha): C_j\to C_i$ for $\alpha: i\to j$, then
the inverse 2-limit $2\pro_{i\in I} C_i$ is the category whose objects are
systems $\{ (E_i)_{i\in I}, (\phi_\alpha)_{\alpha\in {\rm Mor}(I)}\}$
where $E_i\in C_i$ and $\phi_\alpha: C(\alpha)(E_j)\to E_i$ are
isomorphisms, coimpactible in the obvious sense. A morphism in
$2\pro\,  C_i$ from $\{(E_i), (\phi_\alpha)\}$ to 
 $\{(E'_i), (\phi'_\alpha)\}$ is a family $(f_i: E_i\to E'_i)$
compatible with the $\phi_\alpha, \phi'_\alpha$. 
This construction is called quasi-limit in [Gra], p. 217.

\vskip .3cm

\noindent {\bf (A.4) Categories of bundles and torsors.}
Let  $X$ be a finite set. A rank $r$ complex vector bundle on $X$ is
just a collection $(V_x)_{x\in X}$ of $r$-dimensional {\bf C}-vector
spaces. let ${\rm Bun}_r(X)$ be the category of such bundles. 
The correspondence $X\mapsto {\rm Bun}_r(X)$ gives a contravariant
(2-)functor ${\cal S}_0 \to {\cal C}at$, and we extend it to
${\rm Ind}_s^{\aleph_0}({\rm Pro}({\cal S}_0))\simeq {\cal K}$ by taking
the 2-limits, i.e., setting
$${\rm Bun}_r(``\pro"_{i\in I} Y_i) = 
2\ind_{i\in I} {\rm Bun}_r(Y_i), \quad Y_i\in {\cal S}_0, \leqno (A.4.1)$$
$${\rm Bun}_r(``\ind"_{j\in J} Z_j) = 
2\pro_{j\in J} {\rm Bun}_r(Z_j), \quad Z_j\in {\rm Pro}({\cal S}_0), \leqno (A.4.2)$$
Objects of these categories will be called topological {\bf C}-vector bundles.
For a vector bundle $E$ on a finite set $X$ the space of sections
$\Gamma(X,E)$ is an object of ${\rm Vect}_0$, for a bundle $E$ over a profinite
set $Y$ the space $\Gamma(Y,E)$ is an object of 
${\rm Ind}({\rm Vect}_0)={\rm Vect}$ (i.e., just a vector space) and 
for a bundle over an object of $\cal K$ it is an object of
${\rm Pro}_s({\rm Vect})$. The following is elementary.

\proclaim (A.4.3) Proposition. For $Y\in {\rm Pro}({\cal S}_0)$
the category ${\rm Bun}_r(Y)$ is equivalent to the category of
rank $r$ locally constant sheaves of
${\bf C}$-vector spaces on $|Y|\in {\cal P}$. 

Such locally constant sheaves are the same as rank $r$ vector
bundles given by locally constant transition functions. 

The same approach will be followed in all other situations. For example, if
$k$ is a local field such as in (1.2), then for $Y\in {\cal P}$ we have
the category $k-{\rm Bun}_r(Y)$ of topological $k$-vector bundles on $Y$
of rank $r$ (given by continuous, not necessarily locally constant,
transition functions) and their continuous linear morphisms. The category of
topological $k$-vector bundles on $Z=(Z_i)_{i\in I}\in {\rm Ind}_s^{\aleph_0}
({\cal P})$ is defined by
$$k-{\rm Bun}_r(Z) = 2\pro_{i\in I} k-{\rm Bun}_r(Z_i). \leqno (A.4.4)$$

\vskip .1cm

 Let $A$ be a free Abelian group of finite rank. For a finite set $X$ let 
$A-{\rm Tors}(X)$ be the category of $A$-torsors over $X$. We extend
the (2-)functor $A-{\rm Tors}: {\cal S}_0\to{\cal C}at$ to
${\rm Ind}_s^{\aleph_0}({\rm Pro}({\cal S}_0))$ as before:
$$A-{\rm Tors}(``\pro"_{i\in I} Y_i) = 2\ind_{i\in I} A-{\rm Tors}
(Y_i), \quad Y_i\in {\cal S}_0,\leqno (A.4.5)$$
$$A-{\rm Tors}(``\ind"_{j\in J} Z_i) = 2\ind_{j\in J} A-{\rm Tors}
(Z_j), \quad Z_j\in {\rm Pro}({\cal S}_0).\leqno (A.4.6)$$

\proclaim (A.4.7) Proposition. Let $Y=``\pro"_{i\in I} Y_i\in
{\rm Pro}({\cal S}_0)$. For $B=(i, B_i)\in A-{\rm Tors}(Y)$ set
$$|B|= \pro_{(\alpha, j)\in I/i} Y(\alpha)^*B_i,$$
where $I/i$ is the category of arrows $i\buildrel\alpha\over\to j$. Then
$|B|\to |Y|$ is an $A$-torsor in a topological sense. The correspondence
$B\mapsto |B|$ identifies $A{\rm -Tors}(Y)$ with the category of
such topological torsors.

\vskip .2cm

\noindent {\bf (A.5) Function spaces on torsors.} Let $A$ be as before,
$T_A= {\rm Spec} \, {\bf C}[A]$ and $B$ be a principal
homogeneous $A$-space. 
 We denote by
 ${\cal F}(B)$ resp.
${\cal F}_0(B)$ the space of functions, resp. finitely 
supported functions $B\to {\bf C}$. They are modules over ${\bf C}[A]$,
 and ${\cal F}_0(B)$ is free of rank 1.   Let ${\bf C}(T_A)$ be the field of
rational functions on $T_A$. Define
$${\cal F}^{\rm rat}(B) = {\cal F}_0(B)\otimes_{{\bf C}[A]}{\bf C}(T_A).
 \leqno (A.5.1)$$
Let $\Lambda\i A\otimes {\bf R}$ be a strictly convex cone.
Denote by ${\cal F}_\Lambda(B)\i {\cal F}(B)$ the space of functions
whose support is contained in an affine translation of $\Lambda$.
Clearly, ${\cal F}_\Lambda(A)$ is a ring (by convolution)
containing ${\cal F}_0(A)=
{\bf C}[A]$ and ${\cal F}_\Lambda (B)$ is an ${\cal F}_\Lambda(A)$-module.

An element $f: A\to {\bf C}$ of ${\cal F}_\Lambda(A)$ can be seen as a
formal power series $\sum_{a\in A} f(a) t^a$. let ${\cal F}_\Lambda^{\rm rat}
(A)$ be the set of $f$ for which this series is an expansion of a rational
function. The conditions of $f$ imposed by saying that $f\in {\cal F}_\Lambda
^{\rm rat}(A)$, are translationally invariant, so for an $A$-torsor $B$
it makes sense to say that $f\in {\cal F}_\Lambda(B)$ is rational.
 So we get an 
${\cal F}_\Lambda^{\rm rat}(A)$-submodule ${\cal F}_\Lambda^{\rm rat}(B)
\i {\cal F}_\Lambda(B)$. Note that we have a natural embedding
(summation map)
$$\Sigma_\Lambda: {\cal F}_\Lambda^{\rm rat}(B) \hookrightarrow
{\cal F}^{\rm rat}(B).\leqno (A.5.2)$$
Let now $X$ be a finite set and $B$ be an $A$-torsor over $X$. 
 The above constructions, applied pointwise (in $X$) give the
spaces
$${\cal F}_0(B)\in {\rm Mod}^0_{{\bf C}[A]}, \quad {\cal F}(B) \in {\rm Mod}_{{\bf C}[A]}, \quad {\cal F}^{\rm rat}(B)\in {\rm Mod}^0_{{\bf C}(T_A)},$$
$${\cal F}_\Lambda(B) \in {\rm Mod}^0_{{\cal F}_\Lambda(A)}, \quad
{\cal F}_\Lambda^{\rm rat}(B)\in {\rm Mod}^0_{{\cal F}^{\rm rat}_\Lambda
(A)}.\leqno (A.5.3)$$
We then extend these constructions to torsors over objects of 
${\cal K} = {\rm Ind}_s^{\aleph_0}({\rm Pro}({\cal S}_0))$. 
For $Y=(Y_i)_{i\in I}\in {\rm Pro}({\cal S}_0)$ and $B=(i, B_i)\in A-{\rm Tors}(Y)$ set
$$ {\cal F}(B) = ``\ind"_{(\alpha, j)\in I/i} {\cal F}(Y(\alpha)^* B_i) 
\in {\rm Ind}({\rm Mod}_{{\bf C}[A]})
\leqno (A.5.4)$$
and similarly for ${\cal F}_0$, ${\cal F}^{\rm rat}$ etc.
 Note that the
vector space $|{\cal F}(B)|$ is contained in but not, in general,
 equal to, the
space of locally constant functions on the locally compact space $|B|$.
On the other hand, ${\rm Ind}({\rm Mod}^0_{{\bf C}[A]}) =  
{\rm Mod}_{{\bf C}[A]}$ so the object ${\cal F}_0(B)$ of the former category can
be identified with the limit module $|{\cal F}_0(B)|$, and this module coincides with the space with compactly supported locally constant functions on
$|B|$. Invoking (A.1.1)(a), we will consider ${\cal F}_0(B), 
{\cal F}^{\rm rat}(B)$, ${\cal F}_\Lambda(B), {\cal F}_\Lambda^{\rm rat}(B)$
simply as modules over their respective rings. 

Next, let 
$Z=(Z_j)_{j\in J}\in {\rm Ind}_s^{\aleph_0}({\rm Pro}({\cal S}_0))$
and $B\in A-{\rm Tors}(Z)$, so $B$ is a system of torsors $B_i$ over $Z_i$ equipped with identifications $Z(\alpha)^*B_j\to B_i$,
$\alpha: i\to j$. We define the function spaces on $B$ to be the pro-objects
$${\cal F}_0(B) = ``\pro"_{i\in I}{\cal F}_0(B_i) \in {\rm Pro}
({\rm Mod}_{{\bf C}[A]}),\,\,
{\cal F}(B) = ``\pro"_{i\in I} {\cal F}(B_i) \in {\rm Pro}({\rm Ind}
({\rm Mod}_{{\bf C}[A]}))\leqno (A.5.6)$$
and so on. 

\vskip .2cm

\noindent {\bf (A.6) Rational maps of pro-vector bundles.} 
Let $S$ be an irreducible affine algebraic variety over {\bf C}. We denote by
$\Gamma_{\rm reg}: {\rm QCoh}(S)\to {\rm Mod}_{{\bf C}[S]}$ the
 standard equivalence (the functor of regular sections) and 
$\Gamma
_{\rm rat}: {\rm QCoh}(S)\to {\rm Mod}_{{\bf C}(S)}$
functor of rational sections (i.e., the composition of $\Gamma_{\rm reg}$
with the extension of scalars to ${\bf C}(S)$.  For ${\cal F}, {\cal G}\in
{\rm QCoh}(S)$ define
$${\rm Hom}_{\rm rat}({\cal F}, {\cal G}) = {\rm Hom}_{{\bf C}(S)}
(\Gamma_{\rm rat}({\cal F}), \Gamma_{\rm rat}({\cal G})).$$
Suppose ${\cal F}, {\cal G}$ are free (possibly of
infinite rank). Then $R\in {\rm Hom}_{\rm rat}({\cal F}, {\cal G})$
can  given by a (possibly infinite) matrix $\|r_{ij}\|$ over ${\bf C}(S)$.
 For a rational
function $a\in {\bf C}(S)$ we denote by ${\rm Sing}(a)$ its 
singular locus,
i.e., the union of
irreducible hypersurfaces $Z\i S$ such that ${\rm ord}_Z(a) < 0$.
For a rational morphism $R: {\cal F}\to {\cal G}$ of free
${\cal O}_S$-modules and an irreducible
 hypersurface
$Z\i S$ we define ${\rm ord}_Z(R) = \inf_{i,j} {\rm ord}_Z(r_{ij})$,
where $\|r_{ij}\|$ is the matrix of $R$.
This infimum (possibly equal to $(-\infty)$) is independent on the
choice of bases. We set ${\rm Sing}(R)$ to be the union of $Z$ such that
${\rm ord}_Z(R)<0$. 

Let now ${\cal F}= ({\cal F}_i)$, ${\cal G} = 
({\cal G}_j)$ be two strict pro-objects in the category of free
${\cal O}_S$-modules. According to the general definition of
morphisms of pro-objects, we set
$${\rm Hom}_{\rm rat}({\cal F}, {\cal G}) = \pro_j \ind_i {\rm Hom}
_{\rm rat}({\cal F}_i, {\cal G}_j).\leqno (A.6.2)$$
So a rational morphism $A: {\cal F}\to {\cal G}$ is a compactible collection
of rational morphisms $A^{(j)}: {\cal F}_{i(j)} \to
{\cal G}_j$. We define 
$${\rm Sing}(A) = \bigcup_j {\rm Sing}(A^{(j)}, \quad
{\rm ord}_Z(A) = \inf_j {\rm ord}_Z(A^{(j)}).\leqno (A.6.3)$$
It is clear that these concepts are unchanged under left of right
composition with an isomorphism of pro-vector bundles.

\vskip .3cm

\noindent {\bf (A.7) Mellin transform on torsors.} Let $B$ be
a principal homogeneous $A$-space.
For $\lambda\in T_A$ let ${\cal L}(\lambda)_B$ denote the 1-dimensional
space of functions $\psi: B\to {\bf C}$ which are homogeneous of degree
$\lambda$, i.e., satisfy
$$\psi(a+b)=\lambda^a \psi(b), \quad a\in A, b\in B.\leqno (A.7.1)$$
When $\lambda$ varies, the ${\cal L}(\lambda)_B$ form the fibers of an
algebraic line bundle ${\cal L}_B$ on $T_A$. Let $B^- = {\rm Hom}_A(B,A)$
be the torsor dual to $B$ and $i: T_A\to T_A$ be the inversion map:
$i(\lambda) = \lambda^{-1}$. Notice that we have natural isomorphisms of
algebraic line bundles 
$${\cal L}_B^* \simeq {\cal L}_{B^-} \simeq i^*{\cal L}_B.$$
The following is obvious from the construction.

 \proclaim (A.7.2) Proposition. We have a natural identification
$$m_B: {\cal F}_0(B) \to \Gamma_{\rm reg} (T_A, {\cal L}_B^*) 
\simeq \Gamma_{\rm reg}(T_A, i^* {\cal L}_B) = 
\Gamma_{\rm reg}(T_A, {\cal L}_B).
$$ 
 It is defined
as follows. For $\phi\in {\cal F}_0(B)$ the linear form $m_B(\phi):
{\cal O}_B\to {\cal O}_{T_A}$ takes $\phi\in {\cal L}(\lambda)_B$
into $\sum_{b\in B} \phi(b)\psi(b)$. 

 We will call
$m_B$ the {\it Mellin transform} for $B$. Note that $m_B$ induces the
identification
$$m_B^{\rm rat}: {\cal F}^{\rm rat}(B)\to\Gamma_{\rm rat}(T_A, {\cal L}_B)$$

If $X$ is a finite set of, say, $n$ elements, and $B$ is an $A$-torsor
over $X$, then the above considerations can be applied to each fiber
$B_x, x\in X$ and give, for any $\lambda\in T_A$, a line bundle ${\cal L}
(\lambda)\in {\rm Bun}_1(X)$. Note that we have natural
isomorphisms
$${\cal L}(\lambda)\otimes {\cal L}(\lambda') \simeq {\cal L}(\lambda\lambda').\leqno (A.7.3)$$
The space $\Gamma(X, {\cal L}(\lambda)$ (of dimension $n$) will be denoted
$V_\lambda$. When $\lambda$ varies, the $V_\lambda$ form the fibers of an
algebraic vector bundle $V$ on $T_A$ of rank $n$, and we have  natural
identifications
$$m_B: {\cal F}_0(B)\to \Gamma_{\rm reg}(T_A, i^* V) = \Gamma_{\rm reg}(T_A, V), \quad m_B^{\rm rat}: {\cal F}^{\rm rat}(B)\to\Gamma_{\rm rat}(T_A, V).
\leqno (A.7.4)$$
Let $Y=(Y_i)_{i\in I}\in {\rm Pro}({\cal S}_0)$ and $B=(i, B_i)$ be an 
$A$-torsor over $B$.  We get line bundles ${\cal L}(\lambda)$ on $Y$ in the sense of (A.4.1), which satisfy (A.7.3). The space $V_\lambda =
\Gamma(Y, {\cal L}(\lambda))$ is an object of ${\rm Ind}({\rm Vect}_0)$,
so we can consider it as just a vector space $|V_\lambda|$. 
The ``bundle" $V$ is an object
of ${\rm Ind}({\rm Coh}(T_A))$, so we can identify it with the quasicoherent
sheaf $|V|\in {\rm QCoh}(T_A)$. The functor
$$\Gamma_{\rm reg}: {\rm Ind}({\rm Coh}(T_A))\to {\rm Ind}({\rm Mod}^0_{
{\bf C}[A]})$$
obtained by extending $\Gamma_{\rm reg}: {\rm Coh}(T_A)\to {\rm Mod}^0
_{{\bf C}[A]}$ to ind-objects, is identified with the functor of global sections
${\rm QCoh}(T_A)\to {\rm Mod}_{{\bf C}[A]}$ so we have a Mellin
transform $m_B$ as in (A.7.4), with
$\Gamma_{\rm reg}(T_A), V)$ understood in either of two senses. In addition,
taking a trivialization of $B_i$ over $Y_i$, we get the following.

\proclaim (A.7.5) Proposition. For $Y$ as above, the sheaf $|V|$ is free
 so it can be regarded as an algebraic
vector bundle (possibly of infinite rank). The fiber of this bundle at
$\lambda\in T_A$ us canonically identified with the vector space
$|V_\lambda|$.

Next, we extend the above construction to an $A$-torsor $B$ over $Z=(Z_i)_{i
\in I}$ from ${\rm Ind}_s^{\aleph_0} ({\rm Pro}({\cal S}_0))$.  We get line bundles
${\cal L}(\lambda)$ on $Z$ in the sense of (A.4.2), which satisfy
(A.7.3). The ``spaces" $V_\lambda = \Gamma(Z, {\cal L}(\lambda))$
are objects of ${\rm Pro}({\rm Vect})$. The ``bundle" $V$ is an object
of ${\rm Pro}({\rm QCoh}(T_A))$. The functors 
$$\Gamma_{\rm reg}: {\rm Pro}({\rm QCoh}(T_A))\to {\rm Pro}({\rm Mod}_{
{\bf C}[A]}), \quad \Gamma_{\rm rat}: {\rm Pro}({\rm QCoh}(T_A))
\to {\rm Pro}({\rm Mod}_{{\bf C}(T_A)}) \leqno (A.7.6)$$
 are defined componentwise. The
Mellin transforms are now  isomorphisms
$$m_B: {\cal F}_0(B)\to\Gamma_{\rm reg}(T_A, i^*V) = \Gamma_{\rm reg}
(A,V), \quad m_B^{\rm rat}: {\cal F}^{\rm rat}(B)\to \Gamma_{\rm rat}
(T_A, V)
\leqno (A.7.7)$$
in ${\rm Pro}({\rm Mod}_{{\bf C}[A]})$. Unlike the previous case, neither of
the pro-categories involved can be given an ``elementary" description,
so dealing with pro-objects is unavoidable at this stage.

\vfill\eject

\centerline {\bf References}

\vskip 1cm

\noindent [AGV] M. Artin, A.Grothendieck, J.-L. Verdier, SGA4, Theorie
des Topos et Cohomologie Etale des Schemas, t.1, Lecture Notes in Math.
{\bf 269}, Springer-Verlag, 1972.

\vskip .2cm

\noindent [AM] M. Artin, B. Mazur, Etale Homotopy,
Lecture Notes in Math. {\bf 100}, Springer-Verlag, 1969. 

\vskip .2cm

\noindent [B] K. Brown, Buildings, Springer-Verlag, 1989.

\vskip .2cm

\noindent [Cas] W. Casselman, Unramified principal series for p-adic groups I.
The spherical function, 
{\it Compositio Math.} {\bf 40} (1980), 387-406.

\vskip .2cm

\noindent [Ch] I. Cherednik, Double affine Hecke algebras and Macdonald's
conjectures, {\it Ann. Math.} {\bf 141} (1995), 191-216. 

\vskip .2cm

\noindent [CK] N. Chriss, K. Khuri-Makdisi, On the Iwahori-Hecke algebra
of a p-adic group, {\it Int. Math. Res. Notices}, 1998, no.2, 85-100. 

\vskip .2cm

\noindent [Dr] V.G. Drinfeld, Two-dimensional
$l$-adic representations of the fundamental group of a curve over
a finite field and automorphic forms on $GL(2)$, {\it Amer. J. Math.}
{\bf 105} (1983), 85-114.



\vskip .2cm

\noindent [FP] T. Fimmel, A.N. Parshin, Introduction to Higher Adelic
Theory, book in preparation. 

\vskip .2cm

\noindent [GZ] P. Gabriel, M. Zisman, Calculus of Fractions and
Homotopy Theory (Ergebnisse der Math. {\bf 35}), Springer-Verlag, 1967.



 \vskip .2cm

\noindent [GG] H. Garland, I. Grojnowski, Affine Hecke algebras associated
to Kac-Moody groups, preprint q-alg/9508019.

\vskip .2cm

\noindent [GGP] I.M. Gelfand, M.I. Graev, I.I. Piatetski-Shapiro,
Representation theory and automorphic functions, Academic Press, 1969.

\vskip .2cm

\noindent [GKV] V. Ginzburg, M. Kapranov, E. Vasserot, Residue construction
of Hecke algebras, {\it Adv. in Math.} {\bf 128} (1997), 1-19. 



\vskip .2cm

\noindent [Gra] J.W. Gray, Formal Category Theory, Lecture Notes in Math.
{\bf 391}, Springer-Verlag, 1974.



\vskip .2cm

\noindent [Kac] V. Kac, Infinite-dimensional Lie algebras, Camrbidge Univ.
Press, 1984. 

\vskip .2cm

\noindent [Kat] K. Kato, The existence theorem for higher local class
field theory, preprint M/80/43, IHES, 1980. 

\vskip .2cm

\noindent [KL] D. Kazhdan, G. Lusztig, Proof of the Deligne-Langlands
conjecture for Hecke algebras, {\it Invent. Math.} {\bf 87} (1987),
153-215. 



\vskip .2cm

\noindent [Lu 1]
 G. Lusztig,
  Singularities,
character
formula and $q$-analog of weight multiplicity,
 {\it Asterisque} {\bf
101-102} (1983), 208-222.

\vskip .2cm

\noindent [Lu 2]  G. Lusztig, Intersection homology methods
in representation theory, Proc. ICM-90, vol.1, p. 155-174,
Math. Soc. Japan, Tokyo, 1991.

\vskip .2cm

\noindent [Mat] H.Matsumoto, Sur les sous-groupes arithm\'etiques
des groupes semi-simples d\'eploy\'es, {\it Ann. ENS}, {\bf 2} (1969),
1-62. 

\vskip .2cm

\noindent [Mil] J. Milnor, Introduction to Algebraic K-theory,
Princeton Univ. Press, 1972.

\vskip .2cm

\noindent [Pa1] A.N. Parshin, On the arithmetic of 2-dimensional schemes.
I, Repartitions and residues, {\it Russian Math. Izv.} {\bf 40} (1976),
736-773.

\vskip .2cm

\noindent [Pa2]
  A.N. Parshin,  
Vector bundles and arithmetic groups I: The higher Bruhat-Tits tree,
 {\it Proc. Steklov Inst. Math.} {\bf  208} (1995), 212-233,
preprint alg-geom/9605001.

\vskip .2cm

\noindent [PS]
A. Pressley, G. Segal,  Loop Groups,  Clarendon Press, Oxford, 1988.

\vskip .2cm

\noindent [Th] R.W. Thomason, Homotopy limits in the category of
small categories, {\it Math. Proc. Cambridge Phil. Soc.}
{\bf 85} (1979), 91-109.

\vskip 2cm

\noindent {\sl Author's address: Department of Mathematics,
Northwestern University, Evanston IL 60208, email:
\hfill\break
kapranov@math.nwu.edu}

\bye